\documentclass[a4paper]{article}
\usepackage{amsthm,amssymb,amsmath,enumerate,graphicx,psfrag}
\usepackage{subfigure}

\usepackage{graphics}
\usepackage{epstopdf}
\newtheorem{definition}{Definition}

\newtheorem{theorem}[definition]{Theorem}
\newtheorem{thm}[definition]{Theorem}
\newtheorem{corollary}[definition]{Corollary}

\newtheorem{lemma}[definition]{Lemma} \newtheorem{lem}[definition]{Lemma}
\newtheorem{conjecture}[definition]{Conjecture}
 
\newtheorem{claim}[definition]{Claim}

 \newcommand{\emtext}[1]{\text{\em #1}}
\def\weightdeg{{\rm d\bar{e}g}}
\def\avgdeg{{\rm d\bar{e}g}}
\def\SD{~\overline{\!\!S~~}\!\!\!\!D}

 \title{An approximate version\\ of the Loebl-Koml\'os-S\'os conjecture}
 \author{Diana Piguet\footnote{School of Mathematics, University of Birmingham. Most of the work was done while supported by  ITI grant no.~1M00216220808.} \ and Maya Jakobine Stein\footnote{Supported by  FAPESP grant no.~05/54051-9  and Fondecyt grant no.~11090141.}}
  \date{}
 
\begin{document}
 \maketitle
\begin{abstract}
Loebl, Koml\'os, and S\'os conjectured that if at least half of the vertices of
a graph $G$ have degree at least some $k\in\mathbb N$, then every tree with at most~$k$ edges is a subgraph of $G$. Our main result is an
approximate version of this conjecture for large enough $n=|V(G)|$, assumed
that $n=O(k)$.

Our result implies an  asymptotic bound for the Ramsey number of trees. We
prove that $r(\mathcal T_k,\mathcal T_m)\leq k+m+o(k+m)$, as $k+m \rightarrow
\infty$.
\end{abstract}

\section{Introduction}

We explore how certain global assumptions on a graph $G$ ensure the existence of specific subgraphs.
More precisely, we are interested in finding trees as (not necessarily induced)
subgraphs. The main conjecture in our investigations makes, to this end,
assumptions on the median degree of $G$.

\begin{conjecture}[Loebl, Koml\'os, S\'os~\cite{discrepency}]\label{LKS-conj}
Let $k>0$. Then every graph on $n\in\mathbb N$ vertices of which at least $n/2$ have degree
at least $k$, contains as  subgraphs all trees with at most $k$ edges.
\end{conjecture}

The original version  for $k=n/2$ was formulated by Loebl, the generalisation to
arbitrary $k$ is due to Koml\'os and S\'os
(see~\cite{discrepency}). The $n=O(k)$ case of Conjecture~\ref{LKS-conj} is
often referred to as the dense case (otherwise the sparse case).

Our main result is an approximate version of Conjecture~1 for the dense case.

\begin{thm}\label{thm:approx}
For every $\eta,q>0$ there is an $n_0\in\mathbb N$ such that for every
graph $G$ on $n\geq n_0$ vertices and every $k\geq q n$ the
following is true.

If at least ${n}/{2}$ vertices of $G$ have degree at
least $(1+\eta)k$, then $G$ contains all trees with
at most $k$ edges.
\end{thm}

For arbitrary $k$, this has been conjectured
by Ajtai, Koml\'os and Szemer\'edi in~\cite{aks}. They gave a proof for the special case $k=n/2$.

\medskip

The exact version, Conjecture~1, is trivial
for stars, and for trees that consist of two stars with adjacent centres. Bazgan, Li, and Wo\'zniak~\cite{blw}
have proved the conjecture for paths.
The authors of the present paper proved in~\cite{diam5} the Loebl--Koml\'os--S\'os conjecture for trees of diameter at
most $5$.

In Loebl's version with $k=n/2$, the
conjecture has recently been proved by Zhao~\cite{zhao} for large enough
graphs. Extending the methods of Zhao, and of the present paper, the full Loebl--Koml\'os--S\'os conjecture
has been proved very recently for the dense case by 
Hladk\'y together with the first author~\cite{DiaHon},
 and independently, by Cooley~\cite{cooley}.

A generalisation of an example due to Zhao~\cite{zhao} shows that the bound for
the number of vertices of high degree in Conjecture~\ref{LKS-conj} is
asymptotically best possible. 
It cannot be replaced by $n/2-n/(k+1)$, whenever $k+1$ is even and divides $n$ (for bounds in other cases, see~\cite{diam5}).

To see this, construct a graph $G$ on $n$ vertices as follows. Divide $V(G)$ into $2n/(k+1)$ sets $A_i$, $B_i$, so that $|A_i|=(k-1)/2$, and $|B_i|=(k+3)/2$, for $i=1,\ldots ,n/(k+1)$. Insert all edges inside each $A_i$, and insert all edges between each pair $A_i$, $B_i$. Now, consider the tree $T$ we obtain from a star with $(k+1)/2$ edges by subdividing each edge but one. Clearly, $T$ is not a subgraph of $G$.

An interesting folklore observation is the following. Assume that there is a
counterexample to Conjecture~\ref{LKS-conj} for the dense case
that does not contain some tree of order $k+1$. By
taking many copies of $G$, we could then construct a counterexample to
Conjecture~\ref{LKS-conj} for the sparse case.

\medskip

The Ramsey number $r(H,H')$ of two graphs, $H$ and $H'$,
is defined as the minimum integer $n$ such for every graph $G$ of order at least
$n$ either $H$ is a subgraph of~$G$, or $H'$  is a subgraph of the
complement $\bar{G}$ of $G$. Extending this definition, we denote by $r(\mathcal H,\mathcal H')$ the Ramsey number of two classes of graphs, $\mathcal H$ and $\mathcal H'$, that is,
$r(\mathcal H,\mathcal H')$ is the minimum integer $n$ such for every graph $G$ of order at least
$n$ either each graph $H\in \mathcal H$ is a subgraph of~$G$, or each graph $H'\in \mathcal H'$  is a subgraph of the
complement $\bar{G}$ of $G$. We write $r(\mathcal H) $ as shorthand for $r(\mathcal H,\mathcal H)$. 

For $i\in\mathbb N$, let $\mathcal T_i$ denote the class of all
trees of order $i$. Zhao's result implies that the Ramsey number $r(\mathcal
T_{k+1})\leq 2k$, for large $k$. Bounds for Ramsey numbers of trees have been studied for instance in~\cite{treeRamsey}).

In the same way as the bound on $r(\mathcal T_{k+1})$ follows from the Loebl conjecture, one can deduce from Conjecture~1, if true, a bound on $r(\mathcal T_{k+1}, \mathcal T_{m+1})$.
Namely, 
for any colouring of the edges of the complete graph $K_{m+k}$ with two colours,
either half of the vertices have degree $k$ in the subgraph induced by the first colour, or
half of the vertices have degree $m$ in the subgraph induced by the second
colour. So the Loebl--Koml\'os--S\'os conjecture would then imply that
$r(\mathcal T_{k+1},\mathcal T_{m+1})\leq k+m$. This upper bound has been conjectured in~\cite{discrepency}, and it is not difficult to see that the bound is best possible.

 Using Theorem~\ref{thm:approx}, we prove this to be asymptotically true. 

\newcounter{coro}\setcounter{coro}{\value{definition}}
\begin{corollary}\label{rams:number}
 $r(\mathcal T_{k+1},\mathcal T_{m+1})\leq k+m+o(k+m)$, as
 $k+m\rightarrow\infty$.
\end{corollary}

It is not difficult to see that the exact bound of  $r(\mathcal T_{k+1},\mathcal T_{m+1})\leq k+m$ also follows from a positive answer to the  Erd\H os--S\'os conjecture.
This well-known  conjecture states that each graph with average degree greater than~$k-1$ contains all trees
with at most $k$ edges as subgraphs.
For partial results on the Erd\H os--S\'os conjecture, see e.g.~\cite{bradob, sacwoz, woz}.
 Ajtai, Koml\'os, Simonovits and Szemer\'edi proved the Erd\H os--S\'os conjecture for large $n$. (Unfortunately, a manuscript is not available yet.)

\medskip

Our proof of Theorem~\ref{thm:approx} is inspired by the proof of the
approximate version of the Loebl conjecture by Ajtai,
Koml\'os and Szemer\'edi~\cite{aks}. Here also, we use the
regularity lemma followed by a Gallai-Edmonds decomposition of the reduced cluster graph. 
This enables us to find a certain substructure in the cluster graph, which
contains a large matching, and captures the degree condition on~$G$. The tree
is then embedded mainly into the regular pairs corresponding to the matching
edges.

We shall see that in the case that $k\geq  n/2$, it is not difficult to obtain the same structure as in~\cite{aks}. Our proof then follows~\cite{aks}, providing all
details.

In the case that $k< n/2$, however, the situation is more
complex. We will have to content ourselves with a less favourable
structure in the cluster graph, which complicates the embedding of the tree. For a brief
outline of the crucial ideas we then employ, see Section~\ref{sec:ov}. The full
proof is given in the remainder of Section~\ref{sec:proof}.

\medskip

Using similar ideas of proof, we  extend Theorem~\ref{thm:approx} in a different direction.
We pursue the question which other subgraphs are contained
in our graph $G$ from Theorem~\ref{thm:approx}.

 Our third result asserts that we can replace the trees with
bipartite graphs that may have a few more edges than trees.

\begin{theorem}\label{thm:bip}
For every $\eta,q>0$ and for every $c\in \mathbb{N}$ there is an
$n_0\in\mathbb N$ so that for each graph $G$ on $n\geq n_0$ vertices
and each $k\geq qn$ the following is true.

If at least ${n}/{2}$ vertices of $G$ have degree at
least $(1+\eta)k$, then each connected bipartite graph $Q$ on $k+1$ vertices with at most $k+c$ edges is a subgraph of $G$.
\end{theorem}

In particular, the condition of Theorem~\ref{thm:approx} allows for embedding even cycles in~$G$:

\begin{corollary}\label{cor:cyc}
For every $\eta,q>0$ there is an
$n_0\in\mathbb N$ so that for all graphs $G$ on $n\geq n_0$ vertices
and each $k\geq qn$ the following is true.

If at least ${n}/{2}$ vertices of $G$ have degree at
least $(1+\eta)k$, then $G$ contains all even cycles of length at most $k+1$.
\end{corollary}

Theorem~\ref{thm:bip} does not hold for $\eta =0$, as is witnessed by the
following example. Take the complete graph on $k$ vertices and the empty graph on $k$ vertices. Connect these two graphs with a matching of order $k$. The graph we obtain satisfies the condition of the sharp version of Theorem~\ref{thm:bip}, but does not contains the cycle of length $k+1$.

Also, the condition that $Q$ is bipartite is necessary. This can be seen by
considering copies of the complete bipartite graph $K_{(1+\eta)k,(1+\eta)k}$.
This graph satisfies the condition of Theorem~\ref{thm:bip}, but all its subgraphs are bipartite.

\medskip

Our paper is organised as follows. In Section~\ref{sec:reg}, we introduce the regularity lemma and discuss some basic properties of regularity. Our tool for finding the desired
structure of the cluster graph, Lemma~\ref{lem:match}, will be proved in
Section~\ref{sec:match}. All of Section~\ref{sec:proof} is dedicated to the
proof of our main result, Theorem~\ref{thm:approx}.

In Section~\ref{ext}, we explore applications and generalisations of Theorem~\ref{thm:approx}. Our asymptotic bound for Ramsey numbers of trees (Proposition~\ref{rams:number}) will be derived in Section~\ref{sec:ramsey}. In Section~\ref{sec:bip}, we prove
Theorem~\ref{thm:bip}.

\section{Preliminaries}\label{prem}

The purpose  of this section is to introduce the two main tools used in the proofs of Theorem~\ref{thm:approx} and Theorem~\ref{thm:bip}. The first of these tools is the well-known regularity lemma. The second is
Lemma~\ref{lem:match}, which will give structural information on our graph $G$ from Theorem~\ref{thm:approx} (and Theorem~\ref{thm:bip}). We derive it from the Gallai-Edmonds matching theorem.

\subsection{Regularity}\label{sec:reg}

In this subsection, we introduce the notion of regularity, state
Szemer\'edi's regularity lemma, and review a few useful
properties of regularity. All of this is well-known, so
the advanced reader is invited to skip this section.
For an instructive survey on the regularity lemma and its applications, consult~\cite{komlosSimonovits}.

\medskip

Let us first go through some necessary notation.
For a graph $G=(V,E)$, with $W\subseteq E$ and $S\subseteq V$, we will write
$G-W$ for the subgraph $(V,E\setminus W)$ of $G$, and $G-S$ the subgraph of $G$
which is obtained by deleting all vertices of $S$ and all edges incident with vertices of $S$. For
subsets $X$ and $Y$  of the vertex set $V(G)$,  define $N_Y(X)$ as the set of
all neighbours of $X$ in $Y\setminus X$. If $Y=V(G)$, then we omit the index
$Y$ and write $N(X)$. A vertex $x\in V(G)$ is adjacent to the set $Y$ if $xy\in E(G)$ for some $y\in Y$.
If $X$ and $Y$ are disjoint, then let $e(X,Y)$ denote the number of
edges between $X$ and $Y$. The {\em density} of the pair $(X,Y)$ is $d(X,Y):=\frac{e(X,Y)}{|X||Y|}$.

\medskip

A bipartite graph $G$ with partition classes $C_1$ and $C_2$ is called 
$\varepsilon${\em -regular} if for all subsets $C_1'\subseteq C_1$,
$C_2'\subseteq C_2$ with $|C_1'|\geq\varepsilon|C_1|$ and $|C_2'|\geq\varepsilon|C_2|$,
it is true  that $|d(C_1,C_2)-d(C'_1,C'_2)|<\varepsilon$.

A partition $C_0\cup C_1\cup\dots\cup C_N$ of  $V(G)$ is called
{\em$(\varepsilon,N)$-regular}, if
\begin{itemize}
\item $|C_0|\leq\varepsilon n$ and $|C_i|=|C_j|$ for $i, j\in \{1,\ldots,N\}$,
\item all but at most $\varepsilon N^2$ pairs $(C_i,C_j)$ with $i\neq j$ are $\varepsilon$-regular.
\end{itemize}

We are now ready to state Szemer\'edi's regularity lemma.

\begin{theorem}[Regularity lemma, Szemer\'edi~\cite{sze}]\label{lem:reg}
For every $\varepsilon>0$ and $m_0\in\mathbb N$, there exist $M_0,
N_0\in\mathbb N$ so that every graph $G$ of  order $n\geq N_0$ admits an
$(\varepsilon,N)$-regular partition of its vertex set  $V(G)$ with $m_0\leq
N\leq M_0$.
\end{theorem}

Call the partition classes  $C_i$ of $G$ {\em clusters}. Now,
for each graph $G$, for each $(\varepsilon,N)$-regular partition of $V(G)$, and
for any density $p$ define the {\em cluster graph} (sometimes called {\em
reduced graph}) in the following standard way.

First, we construct an auxiliary graph $G_p$ obtained from $G$ by deleting all edges inside the
clusters $C_i$, all edges that are incident with $C_0$,  all edges between
irregular pairs, and all edges between regular pairs $(C_i,C_j)$ of density
$d(C_i,C_j)<p.$ Set  $s:=|C_i|$, and observe that

\begin{equation}\label{eq:G-Gp}
|E(G-G_p)|\leq N\frac{s^2}2+\varepsilon n^2+\varepsilon N^2s^2+\frac{N^2}2
ps^2\leq(\frac{1}{2m_0}+2\varepsilon+\frac{p}{2})n^2.
\end{equation}

Now, the {\em cluster graph} $H=H_p$ on the vertex set $\{ C_i\}_{1\leq
i\leq N}$ has an edge $C_iC_j$ for each pair $(C_i,C_j)$ of clusters that has positive density in $G_p$.
We shall prefer to work with the {\em weighted cluster graph} $\bar H=\bar H_p$ which we
obtain from $H$ by assigning weights $$w(C_iC_j):=d(C_i,C_j) s$$ to the edges $C_iC_j\in E(H)$.

In the setting of weighted graphs,
the (weighted) {\em degree} of a vertex $v$ is defined as $$\weightdeg(v):=\sum_{u\in N(v)}w(vu),$$
and the degree into a subset $U\subseteq V(\bar H)$, where we only count the
weights of edges in $\{v\}\times U$, is denoted by $\weightdeg_U(v)$.
We shall adopt this notation for our weighted cluster graph $\bar H$. For a subset $X\subseteq C_j$, we write
$$\avgdeg_X(C_i):=\frac{e(X,C_i)}{s}=d(X,C)|X|.$$
For a set $\mathcal Y$ of subsets of distinct clusters from $G_p-C_i$, we shall write $\avgdeg_{\mathcal Y}(C_i)$ for $\sum_{Y\in\mathcal Y}\avgdeg_{Y}(C_i)$.

\medskip

We shall often use edges of  $\bar H$ to represent the respective subgraph of 
$G_p$, or sometimes its vertex set. For example, an edge $e=CD\in E(\bar H)$,
might refer to the subgraph of $G_p$ induced by $C\cup D$, or to $C\cup D$ itself. And for a set $U\subseteq C\cup D$, we sometimes use the shorthand $e\cap U$ for $(C\cup D)\cap U$.

\medskip

Let us review some basic properties of $G_p$ and $\bar H$. Let $C,D\in V(\bar H)$:
We call a set $D'\subseteq D$ {\em significant}, if $|D'|\geq \varepsilon s$. A
vertex $v\in C$ is called {\em typical} to a significant set $D'$ if
$\deg_{D'}(v)\geq(d(C,D)-\varepsilon)|D'|$.
Observe that
\begin{equation}\label{eq:for regular pair}
\text{at most }\varepsilon s\text{ vertices of }C\text{ are not typical to a given
significant set }D'.
\end{equation}

Similarly, we have that 
\begin{equation}\label{eq:typicalityFORcluster}
\text{all but at most }\varepsilon s\text{ vertices $v$ of }C\text{ have degree
}\deg_{G_p}(v)\leq \weightdeg(C)+\varepsilon n.
\end{equation}
For proofs of~\eqref{eq:for regular pair} and~\eqref{eq:typicalityFORcluster},
we refer the reader to~\cite{aks}.

Also, almost all vertices
of any cluster $C\in V(\bar H)$ are {\em typical} to almost all significant sets,
in the following sense.
 
If $\mathcal Y$ is a set of significant subsets of clusters in $V(\bar H)$, then
\begin{equation}\label{eq:typicality2}
|\{Y\in \mathcal Y:\; v\text{ is typical to } Y\}|\geq
(1-\sqrt{\varepsilon})|\mathcal Y|,
\end{equation}
for all but at most
$\sqrt{\varepsilon}s$ vertices $v\in C$.

To see this, assume that the set $C'\subseteq C$ of vertices not
satisfying~\eqref{eq:typicality2} is larger than $\sqrt{\varepsilon}s$. Then
\begin{align*}
\sum_{Y\in \mathcal Y}|\{v\in C:\; v\text{ is not typical to }Y\}| \geq &\ \sum_{v\in C'}|\{Y\in \mathcal Y:\; v\text{ is
not typical to }Y\}|  \\
\geq  & \ |C'|\sqrt{\varepsilon}|\mathcal Y|\\  > & \ \varepsilon s |\mathcal Y|.
\end{align*}
Thus there is a $Y\in \mathcal Y$ such that more than $\varepsilon s$
vertices in $C$ are not typical to~$Y$, a contradiction
to~\eqref{eq:for regular pair}.

\subsection{The matching}\label{sec:match}

The main interest in this subsection is Lemma~\ref{lem:match}, which will give us important structural information on the cluster graph $H$ that corresponds to the
 graph~$G$ from Theorem~\ref{thm:approx} (or later Theorem~\ref{thm:bip}). 
 Lemma~\ref{lem:match} appeared in~\cite{aks} but only a weaker variant was proved.

\medskip

For the proof of Lemma~\ref{lem:match}, we need a simplified version of the
Gallai-Edmonds matching theorem, a proof of which can be found for example
in~\cite[p.~41]{diestelBook05}.

A $1$-{\em factor}, or {\em perfect matching}, of a graph $G$ is a $1$-regular spanning subgraph of $G$.
We call $G$ {\em factor-critical}, if for each $v\in V(G)$, there exists a perfect matching of $G-v$.

\begin{thm}[Gallai, Edmonds]\label{lem:galEdm}
Every graph $G$ contains a set $S\subseteq V(G)$ so that each component of $G-S$ is factor-critical, and so that there is a matching in $G$ that matches the vertices of $S$ to vertices of different components of~$G-S$.
\end{thm}

We are now ready for one of the key tools in the proof of
Theorem~\ref{thm:approx}. Recall that we often conveniently use $M$ to
represent $V(M)$.

\begin{lem}\label{lem:match}
Let $\bar H$ be a weighted graph on $N$ vertices, and let $K\in\mathbb R$. Let  $L$ be the set of those vertices $v\in V(\bar H)$ with $\weightdeg(v)\geq K$. If $|L|>N/2$, then there are two
adjacent vertices $v_A,v_B\in L$, and a matching $M$ in $\bar H$ such that one of the following holds.
\begin{enumerate}[(a)]
\item $M$ covers 
$N(\{v_A,v_B\})$,
\item $M$ covers $N(v_A)$, and $\weightdeg_{ M\cup L}(v_B)\geq K/2$. Moreover, each edge in $M$ has at most one endvertex in $N(v_A)$.
\end{enumerate}
\end{lem}

\begin{proof}
Observe that we may assume that $Y:=V(\bar H)\setminus L$ is independent. (In fact, otherwise we simply delete the edges in $E(Y)$, which will not affect the degree of the vertices in $L$.) Now,
Theorem~\ref{lem:galEdm} applied to the unweighted version of $\bar H$ yields
a set $S\subseteq V(\bar H)$. Among all matchings $M'$ satisfying the
conclusion of Theorem~\ref{lem:galEdm} with $S$, choose $M'$ so that it contains a
maximal number of vertices of $Y$.

Set $L':=L\setminus S$. We shall show that either (a) holds or $L'$ is independent. Suppose there is an edge $uv$ with endvertices
$u,v\in L'$. Then $uv$ lies in some component $C$ of $\bar H-S$. If $V(C)\cap V(M')=\emptyset$, let $M''$ be a $1$-factor of $C-u$, and if $V(C)\cap V(M')=\{ x\}$, then let $M''$ be a $1$-factor of $C-x$. In either case~(a) holds for $v_Av_B=uv$ with $M:=M'\cup M''$.
So, from now on, we assume that $L'$ is independent.

Then, each edge of $\bar H$ that is not incident with $S$ has one
endvertex in $L'$, and one in~$Y$. Consider any component $C$ of $\bar H-S$. Since
$C$ is factor-critical, we
have that $|(C-u) \cap Y|=|(C-u)\cap L'|$, for every $u\in V(C)$. Hence, $C$ consists of only one vertex, and so must every component of $\bar H-S$.

\begin{figure}
\centering
\includegraphics[width=0.7\textwidth]{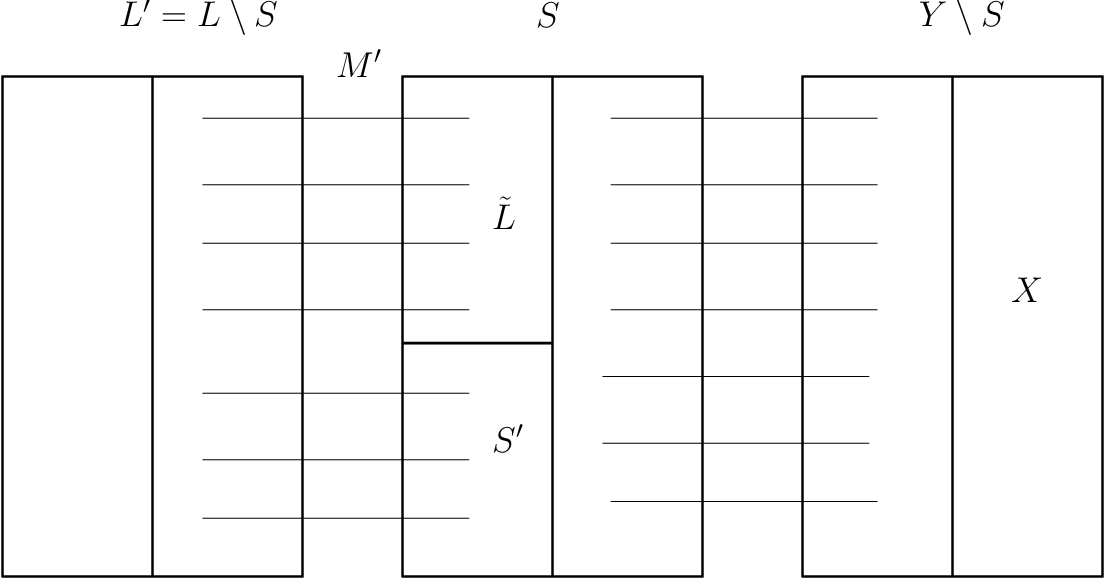} 
\caption{The graph $\bar H$ with the matching $M'$, and sets $L$, $S$ and $Y$.}
\end{figure}


Denote by $X$ the subset of $Y$ that is not covered by $M'$. Set $\tilde
L:=N(L')\cap L\subseteq S$ (see Figure~1).  Now, if there is a vertex $v_B\in\tilde L$ whose weighted
degree into~$\bar H-X$ is at least $ K/2$, then $v_B$, together with any of its
neighbours $v_A$ in $L'$, satisfies~(b) with $M=M'$. So, we may assume that for each $u\in\tilde L$,
\begin{equation}\label{eq:tildeL}
\weightdeg_{\bar H-X}(u)<K/2,
\end{equation}
and hence $\weightdeg_X(u)\geq K/2$.

On the other hand, $\weightdeg_{\tilde L}(w)<K$ for each $w\in X$. Thus, by double (weighted) edge-counting, it follows that
\begin{equation}\label{eq:tilde}
|X|\geq\frac{|\tilde L|}2.
\end{equation}

Set $S':= S\cap Y$.
By~\eqref{eq:tildeL}, the total weight of the edges in $E(\tilde L\cup S',L')$ is less than $|\tilde L| K/2+|S'|K$, while each vertex of $L'$ has weighted degree at least $K$ into $\tilde L \cup S'$. Thus, again by double edge-counting, and by~\eqref{eq:tilde},
\begin{equation}\label{nine}
|X|+|S'|\geq \frac{|\tilde L|} 2+|S'|>|L'|.
\end{equation}

Furthermore, since~$Y$ is independent, $M'$ matches $S'\subseteq Y$ to $L'$.
Thus $|L'|\geq |S'|+|L\setminus M'|$, and so, by~\eqref{nine},
\[
|X|>|L\setminus M'|.
\]

Since $|L|>\frac N2$, this implies that $M'$ contains an edge $uv$ with both $u,v\in L$. We may assume that $v\in L'$ and $u\in \tilde L$. By~\eqref{eq:tildeL}, $u$ has a neighbour~$w$ in $X$. Hence, the matching $M'\cup \{uw\}\setminus \{uv\}$ covers more vertices of $Y$ than $M'$ does,
a contradiction to the choice of~$M'$.
\end{proof}

Note that in the case $K\geq N/2$ the situation in Lemma~\ref{lem:match} is less complicated. In that case, observe that clearly $|S|\leq |V(\bar H-S)|$. So, either $|S|=|V(\bar H-S)|$ (in which case  conclusion (a) of Lemma~\ref{lem:match} holds), 
or there is a component $C$ of $\bar H-S$ that has more than one vertex. Thus,
as $C$ is factor-critical, there exists an edge in $C\cap (L'\times L')$, and
(a) holds again.
%
In the case $k\geq n/2$, this observation simplifies our proof of 
Theorem~\ref{thm:approx} considerably, as then only the simplest case needs to
be treated.

\section{Proof of Theorem~\ref{thm:approx}}\label{sec:proof}

The organisation of this section is as follows. The first subsection is devoted
to an outline of our proof, highlighting the main ideas, leaving out all
details. In Subsection~\ref{prep}, assuming that we are given a host graph $G$
and a tree $T^*$ as in Theorem~\ref{thm:approx}, we shall first apply the
regularity lemma to $G$. We then use Lemma~\ref{lem:match} to find a
suitable matching of the corresponding weighted cluster graph~$\bar
H$, which will  facilitate the embedding of $T^\ast$.

We shall prepare $T^\ast$ for this by cutting it into small pieces in
Subsections~\ref{part} and~\ref{switch}. Then, in Subsection~\ref{num-part},
we partition the matching given by Lemma~\ref{lem:match}, according to the
decomposition of the tree $T^\ast$. In Subsection~\ref{sec-embeddingLemma}, we
expose tools that we need for our embedding. What remains is the actual embedding
procedure, which we divide into the two cases given by Lemma~\ref{lem:match}, and treat separately in Subsections~\ref{sec:1} and~\ref{sec:2}.

\subsection{Overview}\label{sec:ov}

In this subsection, we shall give an outline of our proof of Theorem~\ref{thm:approx}.
So, assume that we are given $\eta>0$ and $q>0$. The regularity lemma applied to parameters depending on $\eta$ and $q$ yields an $n_0\in \mathbb{N}$.
Now, let $n\geq n_0$, let $k\geq qn$, let $G$ be a graph of order $n$ that satisfies the condition of Theorem~\ref{thm:approx}, and let~$T^\ast$ be a tree with $k$ edges. We wish to find a subgraph of $G$ that is isomorphic to~$T^\ast$, i.e. we would like to {\it embed} $T^\ast$ in $G$.

In order to do so, consider the weighted cluster graph
$\bar H$ corresponding to $G$ that is given by the regularity lemma.
Denote by $L\subseteq V(\bar H)$ the set of those clusters that have degree at least $(1+\pi')k$ in~$\bar H$, where $\pi'=\pi'(\eta,q)>0$. 
The weighted cluster graph $\bar H$ inherits properties from $G$ resulting in
the fact that $|L|>|V(\bar H)|/2$. Apply Lemma~\ref{lem:match} to $\bar H$ and
$K:=(1+\pi')k$ which yields vertices $A,B\in V(\bar H)$ and a matching $M$. The rest of our proof will be divided into two cases, corresponding to the two possible conclusions  (a) and (b) of Lemma~\ref{lem:match}.


\medskip

If the output of Lemma~\ref{lem:match} is Case~(a), then we shall decompose $T^*$ into small subtrees (of order much below~$\eta k$) and a small set~$SD$ of vertices (of constant order in $n$), so that between any two of our subtrees
lies a vertex from $SD$ (the name $SD$  stands for `seeds'). In fact, $SD$ is the disjoint union of two sets
$SD^A$ and $SD^B$, and each tree $T$ of $T^\ast-SD$ is adjacent to
only one of these two sets, that is, either $N(SD^A)\cap V(T)=\emptyset$ or $N(SD^B)\cap
V(T)=\emptyset$. Denote the set of trees adjacent to $SD^A$ by
$\mathcal {T}_A$, and the set of trees adjacent to $SD^B$ by~$\mathcal {T}_B$.
The formal definition of $SD$, $\mathcal T_A$ and $\mathcal T_B$ can be found in Section~\ref{part}.

Next, in Section~\ref{num-part}, we partition the matching~$M$ from Lemma~\ref{lem:match} into~$M_A$ and~$M_B$. This is done in a way so that  $\avgdeg_{M_A}(A)$ is large enough so that
$F_A:=\bigcup\mathcal {T}_A$ fits into $M_A$, and $\avgdeg_{M_B}(B)$ is large enough  so that $F_B:=\bigcup\mathcal {T}_B$ fits into $M_B$. 

Finally, in Section~\ref{sec:1}, we embed $SD^A$ in $A$ and $SD^B$ in $ B$ and use the regularity of the
edges in $\bar H$ to embed the small trees of $\mathcal {T}_A\cup\mathcal
{T}_B$, one after the other, levelwise, into $M_A\cup M_B$.
The order of this embedding procedure will be such that the already embedded part of~$T^\ast$ is always connected.

Moreover,
the structure of our decomposition of $T^\ast$, and the fact that we embed the
trees from $\mathcal T_A\cup\mathcal T_B$ in the matching edges, ensures that
the predecessor of any vertex $r\in SD^A \cup SD^B$ is embedded in a cluster
that is adjacent to $A$, respectively to $B$ (in which we wish to embed $r$).
This enables us to embed all of~$SD$ in $A\cup B$, as planned.

An important detail of our embedding technique is that we shall always try to
{\em balance} the embedding in the matching edges, in the sense that the used
part of either endcluster should have about the same size. 
We only allow for an unbalanced embedding if the degree of $A$ resp.~$B$ into
one of the endclusters of the concerned edge is already `exhausted' (cf.\
Property~$(\diamond)$ in Section~\ref{sec-embeddingLemma}). In practice, this means that
whenever we have the choice into which endcluster of an edge $e\in M$ we embed the root of some tree of $ \mathcal T_A\cup  \mathcal T_B$, we shall choose the side carefully.

In this manner, we can ensure that all of $T^\ast$ will fit into $M$ (or more precisely into the corresponding subgraph of $G$). This finishes the embedding of $T^\ast$ in Case~(a) of Lemma~\ref{lem:match}.

\medskip

In Case~(b) of Lemma~\ref{lem:match},  it is not possible to partition the matching $M$ into~$M_A$ and~$M_B$ so that  $F_A$ fits into $M_A$ and $F_B$ fits into $M_B$, as in Case~(a). More precisely, for any partition of $M$ into $M_A$ and $M_B$, if $\avgdeg_{M_A}(A)$ allows for the embedding of a forest of order $t$, say, in $M_A$, then $\avgdeg_{M_B\cup L}(B)$ only guarantees for the embedding of a forest of order at most $(k-t)/2$ in the subgraph of $G_p$ induced by $M_B$ and the edges incident with~$L'$, where $L':=L\setminus M$. For more details on this, see Lemma~\ref{lem:numbers}.

We use a combination of two strategies to overcome this problem.
Firstly, we shall embed~$T^\ast$
in two phases, leaving for the second phase some subtrees that are (each)
adjacent to only one vertex from $SD$. Secondly, we shall embed some of the
trees from $\mathcal T_B$ in part of the matching reserved for $F_A$. This
means that we `switch' some of our trees to $\mathcal T_A$.

\smallskip

Let us explain the two strategies in more detail. We modify our sets $\mathcal
{T}_A, \mathcal {T}_B$, in the following way. Denote by $\bar{\mathcal  T}_A$
the set of those trees from $\mathcal {T}_A$ that are adjacent to only one
vertex from~$SD^A$, and similarly define $\bar{\mathcal  T}_B$. (Observe
that~$T^*$ remains connected after deleting any tree in $\bar{\mathcal T}_A\cup\bar{\mathcal T}_B$.)

We may assume that $$|V(\bigcup\bar{\mathcal  T}_A)|\geq
|V(\bigcup\bar{\mathcal  T}_B)|.$$ Finally, set $\mathcal {T}':=(\mathcal
{T}_A\cup \mathcal {T}_B)\setminus (\bar{\mathcal  T}_A\cup \bar{\mathcal 
T}_B)$. Our plan now is to first embed the trees from $\mathcal  T'\cup\bar
{\mathcal T}_B$ together with the vertices from $SD$ and to postpone the
embedding of $\bar F_A:=\bigcup\bar{\mathcal  T}_A$ to a later stage. 
As the part of the tree embedded in the first phase is connected,
 we avoid the difficulty of having to connect already embedded parts of $T^\ast$  in the second phase.

Now, we shall partition~$M$ into~$M_f$ and~$\bar M_B$ so that
$\avgdeg_{M_f}(A)$ allows for the embedding of  $\bigcup\mathcal {T}'$, and
$\avgdeg_{\bar M_B\cup L}(B)$ allows for the embedding of $\bar F_B:=\bigcup
\bar{\mathcal  T}_B$. This actually means that the place we reserved for the
embedding of $F_B- \bar{F}_B$ lies in $M_f$. Therefore, we shall `switch' this
forest to $\mathcal {T}_A$ (which is the second of our strategies).

Let us explain what we mean by \emph{switching}. For each tree $T\in \mathcal {T}_B\setminus \bar{\mathcal  T}_B$, delete all vertices from $T$ that are adjacent to $SD^B$ in $T^\ast$ and add them to $SD^A$. Put the components of what remains of  $T$ into $\mathcal {T}_A$.
Denote the thus enlarged $SD^A$ by $\SD^A$ and set $\SD:= \SD^A\cup
SD^B$.

After switching all trees $T\in \mathcal {T}_B\setminus \bar{\mathcal  T}_B$, denote by $\mathcal {T}_f$ the (enlarged)  set  $\mathcal {T}_A\setminus\bar{\mathcal  T}_A$. That is, $\mathcal  T_f$ consists of all trees from the original $\mathcal  T_A\setminus  \bar{\mathcal  T}_A$, together with all trees we generated by switching. 
It will be easy to verify that the switching procedure does not increase too much the number of seeds.

Also, each tree from $\mathcal {T}_f$ and $\bar{\mathcal  T}_A$ is adjacent
only to the enlarged $\SD^A$, and each tree from  $\bar{\mathcal  T}_B$ is still adjacent only to $SD^B$. For details on the switching procedure, consult Section~\ref{switch}.

It remains to embed $T^*$ in $G$, which is done in Section~\ref{sec:2}.
We first embed the vertices from $\SD^A\cup SD^B$ in $A\cup B$, embed
$F_f:=\bigcup\mathcal {T}_f$ in~$M_f$, and embed part of $\bar {\mathcal T}_B$
in $\bar M_B$, in the same way as in Case~(a). In a second phase, we embed the
remaining trees from $\bar{\mathcal T_B}$  into edges of $H$ that are incident
with $L'$. For each tree, we are able to find a free space in a
suitable edge because of the high degree of the clusters  from $L'$. 

\medskip

In the remaining third phase we wish to embed $\bar F_A$. We shall now use all
of~$M$, forgetting about the partition into $M_f$ and $\bar M_B$. The neighbours of the trees from $
\bar{\mathcal  T}_A$ in~$\SD^A$ have already been embedded in the first
phase. Having chosen their images carefully then, ensures that now they have
still large enough degree into what is not yet used of $M$. Hence, there is
enough place for $\bar F_A$ in $M$.

Also,  it is essential here that each edge of $M$ meets~$N(A)$ in at most one
cluster. The reason is that parts of these clusters might have been used in
the first and second phases of the embedding. So, some of the edges involved
might be unbalanced, in the sense above, because e.\,g.~the degree of $B$ was
such that we were not able to choose the endcluster in which we embedded the
roots of the trees from $ \bar{\mathcal T}_B$. 
However,  as each edge of~$M$ has at most one endcluster in~$N(A)$,
 it is irrelevant whether the embedding is balanced or not in these edges.

The embedding itself of $\bar F_A$ is done as before. This finishes the sketch of our proof in Case~(b).

\subsection{Preparations}\label{prep}

We shall now prove Theorem~\ref{thm:approx}.
First of all, we fix a few constants depending on $\eta$ and $q$.
Set $$\pi:=\min\{\eta,q\}, \quad \varepsilon:=\frac{\pi^7q}{25\cdot10^7}\quad \text{ and }\quad m_0:=\frac
{500}{q\pi^3}.$$

The regularity lemma (Theorem~\ref{lem:reg}) applied to $\varepsilon$,
 and $m_0$ yields natural numbers~$M_0$ and $N_0$.
Fix $$\beta:=\frac{\varepsilon}{M_0},\quad p:=\frac{\pi^3q}{250}\quad \text{ and }\ \ \
n_0:=\max\left\{N_0,\frac {64 M_0}{\beta p}
\right\}.$$

Thus our constants satisfy the following relations

\begin{equation*}\label{1}
\frac 1{n_0}\ \ll\ \beta\  \ll\  \varepsilon\  \ll\ \frac{1}{m_0}\ <\  p\ \ll\ \pi\ \leq\  q,
\end{equation*}

where $a\ \ll \ b$ stands for the fact that $a<\frac {\pi }{100}b$.

In particular, $p$ satisfies
\begin{equation}\label{eq:p}
4\varepsilon+\frac{1}{m_0}<p.
\end{equation}

Let $n\geq n_0$, let $k\geq qn$, and let $G$ be a graph of order $n$ which has
at least $\frac n2$ vertices of degree at least $(1+\eta)k$. Suppose $T^\ast$ is a tree of order~$ k+1$.
Our aim is to find an embedding $\varphi:V(T^\ast)\to V(G)$ that preserves adjacency.

Now, by Theorem~\ref{lem:reg} there exists an $(\varepsilon,N)$-regular
partition of $V(G)$, with $m_0\leq N\leq M_0$. As in Section~\ref{sec:reg}, let $G_p$ be the subgraph of $G$ that preserves exactly the edges between regular pairs of density at least $p$.

By~\eqref{eq:G-Gp} and by~\eqref{eq:p},
\[
|E(G-G_p)|<pn^2\leq\frac{\pi^3 }{250}kn.
\]
Thus, for all but at most $\frac{\pi^2 }{50} n$ vertices $v$, we have that $\deg_{G_p}(v)\geq \deg_G(v)-\frac{\pi }5 k$. Hence,
\[
G_p\emtext{ has at least } \left (1-\frac {\pi^2}{25}\right)\frac n2\emtext{
vertices of degree at least  }\left (1+\frac{4\pi}{5}\right)k.
\]

Let $\bar H=\bar H_p$ be the weighted cluster graph corresponding to $G_p$. Denote by $L$ the set of those  clusters in $V(\bar H)$ that contain more
than~$\varepsilon s$ vertices of degree at least $(1+\frac{4\pi}{5})k$ in
$G_p$. A
simple calculation shows that $|L|>(1-\frac{\pi^2}{5})\frac{N}{2}$.

 Now, delete
$\min\{\pi^2 N/5, |V(\bar H)\setminus L|\}$ clusters in $V(\bar
H)\setminus L$ to obtain a subgraph of the cluster graph $\bar H$. 
As this subgraph is very similar (or identical) to $\bar H$, in the rest of the
text we shall denote it as well by $\bar H$. So from now on,
by $\bar H$, we shall always refer  to this subgraph. Each vertex in $\bigcup
L$
 drops its degree by at most $\frac {\pi^2}5Ns\leq \frac {\pi k}5$.
Thus, by~\eqref{eq:typicalityFORcluster}, each  $X\in L$ has degree
\begin{equation}\label{Ldeg}
\weightdeg_{\bar H}(X)\ge (1+\frac{3\pi}{5})k-\varepsilon n>
(1+\frac{\pi}{5})k.
\end{equation}
Then Lemma~\ref{lem:match} applied to $\bar H$ and $K:=(1+\frac{\pi}{5})k$ yields an edge $AB\in E(\bar H)$ with $A,B\in L$,  together with a matching $M'$ of $\bar H$, which satisfy~(a) or~(b) of Lemma~\ref{lem:match}.
Obtain $M$ from $M'$ by deleting all edges from $M'$ that are incident with $A$ or
with $B$. If $AA',BB'\in M'$, then $M$ misses
$A$, $A'$, $B$, and $B'$, thus at most three clusters from $N(A)$,
resp.\ from $N(B)$. 
In Case~(a) of Lemma~\ref{lem:match}, we calculate that
\begin{align}\label{eq:degCaseA}
\min\{\weightdeg_{M}(A),\weightdeg_{M}(B)\} &
\geq\  (1+\frac{\pi}{5})k-\frac {3n}N \notag \\
&\geq \ (1+\frac{\pi}{5}-\frac 3{qm_0})k\notag \\
&\geq \ (1+\frac{\pi}{10})k.
\end{align}
Similarly, in Case~(b) it follows that
\begin{equation}\label{eq:degCaseB}
\weightdeg_{M}(A)\geq (1+\frac{\pi}{10})k\ \ \text{ and }\ \
\weightdeg_{M\cup L}(B)\geq (1+\frac{\pi}{10})\frac k2.
\end{equation}

Thus, for the remainder of our proof of Theorem~\ref{thm:approx} we shall work with the assumption that there is a matching $M$ of $\bar H$ and vertices $A,B\notin V(M)$ so that
\begin{enumerate}[1.]
\item\label{case1'}$\avgdeg_M(A),\avgdeg_M(B)\geq (1+\frac{\pi}{10})k$, or
\item\label{case2'}$\avgdeg_M(A)\geq (1+\frac{\pi}{10})k$, $\avgdeg_{M\cup L}(B)\geq (1+\frac{\pi}{10})\frac{k}{2}$, and each cluster in $N(A)$ meets a different edge of $M$.
\end{enumerate}

We shall refer to these two cases as `Case 1' and `Case 2', respectively.
We will embed the tree $T^*$ in the subgraph of $G_p$
corresponding to $\bar H$, using two different strategies in Case~1 and in Case~2.

\subsection{Partitioning the tree}\label{part}

In this section, we shall cut our tree into small pieces. More precisely, we shall define a set $SD\subseteq V(T^\ast)$, and sets $\mathcal T_A$ and $\mathcal T_B$ of disjoint small subtrees of $T^\ast$ which are connected through the vertices from $SD$. Moreover, $SD$ together with the union of all trees from $\mathcal T_A\cup \mathcal T_B$ will span $T^\ast$.

\medskip

Fix a vertex $R$ of $T^\ast$ as the root and regard $T^*$ as a poset having $R$
as the minimal element. For a vertex $x$ of a subtree $T\subseteq T^\ast$,
denote by $T(x)$ the subtree of $T$ induced by $x$ and all vertices $y$ greater than $x$ in the tree-order of~$T^\ast$. (That is, $T(x)$ contains all vertices $y$ such that the path between the root~$R$ and $y$ contains the vertex $x$.) If $R\notin V(T)$, then define  the {\em seed} $sd(T)$ of $T$ as the  maximal vertex of $T^\ast$ which is smaller than every vertex of~$T$.

Our sets $SD=SD^A\cup SD^B$, $\mathcal T_A$ and $\mathcal T_B$  will satisfy:
\begin{enumerate}[(I)]
\item $SD^A\cap SD^B=\emptyset$,
\item  $R\in SD^ A$,  and  $r\in SD$ lies at even distance to $R$ if and only if $r\in SD^A$,
\item $\mathcal T_A\cup \mathcal T_B$ consists of the components of $T^\ast-SD$,
\item $|V(T)|\leq \beta k$, and $sd(T)\in SD$, for each  $T\in \mathcal T_A\cup \mathcal T_B$,
\item\label{V} $\max\{|SD^A|,|SD^B|\}\leq\frac 2\beta$, and
\item  \label{VI}$e_{T^*}(V(F_A),SD^B)=0$, and $e_{T^*}(V(F_B),SD^A)=0$,
\end{enumerate}
where $F_A:=\bigcup_{T\in \mathcal T_A}T$ and $F_B:=\bigcup_{T\in \mathcal T_B}T$ are the forests spanned by $\mathcal T_A$ and~$\mathcal T_B$.

Let us first define $SD$. To this end, we shall inductively find vertices $x_i$, and define auxiliary trees $T^i\subseteq T^\ast$. Set $T^0:=T^\ast$.

\begin{figure}[htp]
     \centering
     \subfigure[Suppose that $x_1,x_4,x_5,x_7,x_9$ are in $T^{i-1}(x_i)$.]{
          \label{fig:b1}
          \includegraphics[width=.44\textwidth]{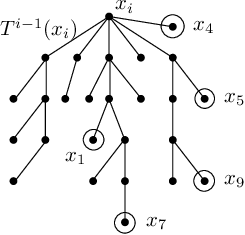}}
     \hspace{.08in}
     \subfigure[Say $x_i\in A'$. Then
     $x_5,x_9\in A'$ and $x_1,x_4,x_7\in B'$, which we mark by circles and squares respectively.]{
          \label{fig:b2}
          \includegraphics[width=.44\textwidth]{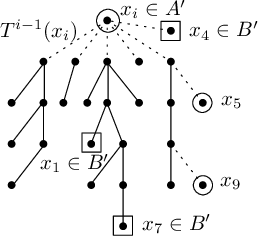}}
     \hspace{.02in}
     \subfigure[We add $y$ and $z$ to $A(T)$. Then $T^{i-1}(x_i)-SD\subseteq \mathcal T_A$. ]{
           \label{fig:b3}
           \includegraphics[width=.44\textwidth]{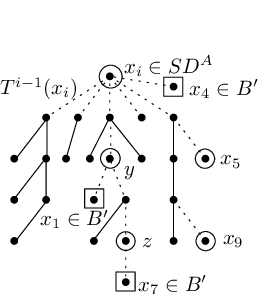}}
     \caption{Phases of the partition of $T^\ast$.}
     \label{fig:2858multifig}
\end{figure}

In step $i\geq 1$, let $x_i\in V(T^*)$ be a maximal vertex in the tree-order of $V(T^{i-1})$ with
\begin{equation}\label{eq-betakTree}
|V(T^{i-1}(x_i))|> \beta k,
\end{equation}
as illustrated in Figure~\ref{fig:b1},
and define $$T^i:=T^{i-1}-(T^{i-1}(x_i)-x_i).$$

Hence,  
\begin{equation}\label{betak-1}
|V(T^{i-1})|-|V(T^i)|>(\beta k-1).
\end{equation}

If there is no  vertex satisfying~\eqref{eq-betakTree}, then set $x_i:=R$, and stop the definition process.
Say our process stops in some step $j$. Let $A'$ be the set of all $x_i$, $i\leq j$, with even distance to the root $R$, and let $B'$ be the set of all other $x_i$. 

Then, by~\eqref{betak-1} and by the definition of $n_0$,
\[
j-1\leq\frac {|V(T^*)|}{\beta k -1}= \frac{k+1}{\beta k -1}\leq \frac
{3}{2\beta}.
\]
 
Hence,
\begin{equation}\label{eq:2/beta}
|A'\cup B'|\leq \frac 2\beta.
\end{equation}

For the sake of condition~(VI), we shall now add a few more vertices to our sets~$A'$ and $B'$, which will result in the desired $SD$.

Let $\mathcal{C}$ be the set of the components of $T^\ast-(A'\cup B')$.
For each  $T\in\mathcal C$ with $sd(T)\in A'$, denote by $A(T)$ the set of vertices  of $T$ that are adjacent to $B'$. Similarly, if $sd(T)\in B'$, then denote by $B(T)$ the set of  vertices of $T$ that are adjacent to $A'$ (cf.~Figure~\ref{fig:b2}). Set $$SD^A:=A'\cup \bigcup_{T\in \mathcal{C}}A(T),\  \text{ and }\ SD^B:=B'\cup\bigcup_{T\in\mathcal{C}}B(T)$$ and set $SD:=SD^A\cup SD^B.$

Since each vertex in $B'$ has at most one neighbour in the union of the $A(T)$, it follows that 
\[
|SD^A\setminus A'|\leq |B'|,
\]
and analogously,
\[
|SD^B\setminus B'|\leq |A'|.
\]
Thus, 
\begin{equation}\label{eq:SD}
\max\{|SD^A|,|SD^B|\}\leq |A'\cup B'|.
\end{equation}

Finally, we shall define $\mathcal T_A$ and $\mathcal T_B$.  Let $\mathcal C'$ be the set of the components of $T^\ast-SD$. Set
\[
\mathcal T_A:=\{T\in \mathcal{C'}:\ sd(T)\in SD^A\}\ \ \ \text{ and }\ \ \ \mathcal T_B:=\{T\in \mathcal{C'}:\ sd(T)\in SD^B\},
\]
as shown in~Figure~\ref{fig:b3}, and define the forests 
\[
F_A:=\bigcup_{T\in \mathcal T_A}T\ \ \ \text{ and }\ \ \ F_B:=\bigcup_{T\in \mathcal T_A}T.
\]

Observe that Conditions (I)--(IV) and (VI) are clearly met and that~\eqref{V} holds because of~\eqref{eq:2/beta} and~\eqref{eq:SD}.
 
This finishes our manipulation of the tree $T^\ast$ in Case~1.

\subsection{The switching}\label{switch}

In Case~2 from Section~\ref{prep}, we shall not only cut our tree to small pieces (cf. Section~\ref{part}), but also switch
some of our small subtrees from one of the two sets $\mathcal T_A$, $\mathcal T_B$ to the other.
 We achieve this by adding some more vertices to $SD$, thus naturally refining our partition of $T^\ast$.

\medskip

 Set
\begin{align*}
\bar{\mathcal T}_A:= &\{T\in \mathcal T_A: e(V(T),SD-sd(T))=0\},\ \text{ and}\\
\bar{\mathcal T}_B:= &\{T\in \mathcal T_B: e(V(T), SD-sd(T))=0 \}.
\end{align*}

We may assume that
\begin{equation}\label{barA}
|\bigcup_{T\in\bar{\mathcal T}_A}V(T)\ |\geq |\bigcup_{T\in\bar{\mathcal T}_B}V(T)\ |.
\end{equation}

\begin{figure}[htp]
     \centering
     \subfigure[A tree $T\in \mathcal T_B\setminus \bar{\mathcal T}_B$, with $sd(T),x_1,x_2,x_3,x_4\in SD^B$.]{
          \label{fig:s1}
          \includegraphics[width=.42\textwidth]{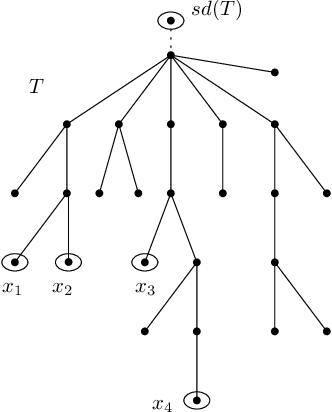}}
     \hspace{.4in}
     \subfigure[The set $S(T)=\{y_1,\ldots,y_4\}$, and the subtrees of $T$ generated by the switching.]{
          \label{fig:s2}
          \includegraphics[width=.42\textwidth]{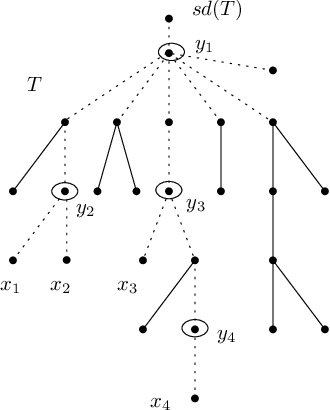}}
     \caption{The switching procedure.}
     \label{fig:switch}
\end{figure}

Now, consider a tree $T\in \mathcal T_B\setminus\bar{\mathcal T}_B$ as in
Figure~\ref{fig:s1}. By~\eqref{VI}, no vertex in $V(T)$ is adjacent to any
vertex in $SD^A$ in $T^*$. Denote by $S(T)$ the set of all vertices in $V(T)$
that in $T^\ast$ are adjacent  to some vertex of  $SD^B$. For illustration see Figure~\ref{fig:s2}.
Set $$\SD^A:=SD^A\cup \bigcup_{T\in \mathcal T_B\setminus \bar{\mathcal T}_B}
S(T) \ \ \ \text{ and }\ \ \SD:=\SD^A\cup SD^B.$$

Finally, define
\[
\mathcal T'_A:=\bigcup_{T\in \mathcal T_B\setminus\bar{\mathcal T}_B}\{C:C\text{ is a component of }T-S(T)\} 
\]
and
\[
\mathcal T_f:=(\mathcal T_A\setminus \bar{\mathcal T}_A)\cup \mathcal T'_A.
\]
(The $f$ in $\mathcal T_f$ stands for `first', as this part of the tree is to be embedded first.)
Finally, set 
\begin{align*}
F_f:= &\bigcup_{T\in \mathcal T_f}T, \\
\bar F_A:=\bigcup_{T\in \bar{\mathcal T}_A}T\ \ & \text{ and }\ \ \bar F_B:=\bigcup_{T\in \bar{\mathcal T}_B}T.
\end{align*}

Observe that our sets $\SD=\SD^A\cup SD^B$, $\mathcal T_f\cup
\bar{\mathcal T}_A$, and $\bar{\mathcal T}_B$ still satisfy conditions (I)-(IV)
and (VI) from Section~\ref{part} (with $SD$, $SD^A$, $\mathcal T_A$, $\mathcal
T_B$, $F_A$, and $F_B$ replaced by $\SD$, $\SD^A$, $\mathcal T_f\cup \bar{\mathcal T}_A$, $\bar{\mathcal T}_B$, $\bar F_A,$ and $\bar F_B$, respectively). Instead of (V), we now have the similar

\begin{enumerate}[(I)']\setcounter{enumi}{4}
\item $|\SD|\leq\frac8\beta$,
\end{enumerate}
since by the definition of $\SD^A$ we know that for each vertex $x$ of $SD^B$,
we have created at most $2$ vertices of $\SD^A\setminus SD^A$ (between $x$ and
the next vertex of $SD^B$ in direction of the root $R$). Thus,
 \[
 |\SD^A|\leq |SD^A|+2|SD^B|\leq \frac 6\beta,
\]
as needed for (V)'.

\subsection{Partitioning the matching}\label{num-part}

In this subsection, we shall divide the matching $M$ into two parts, into which we will later embed the two forests $F_A$, $F_B$, respectively $F_f$ and $\bar F_B$, of $T^*$ that we defined in Subsection~\ref{part}, resp.~in Subsection~\ref{switch}. (The forest $\bar F_A$ will be embedded later).

\medskip

For this, we will need the following number-theoretic lemma, which appeared
also in~\cite{aks}. We give a short proof.

\begin{lem}\label{lem:numbers}
Let $I$ be a finite set, and  let $a,b,\Delta>0$. For $i\in I$, let $a_i,b_i\in
(0,\Delta]$. Suppose that
\begin{equation}\label{sum}
\frac{a}{\sum_{i\in I}a_i}+\frac{b}{\sum_{i\in I}b_i}\leq 1.
\end{equation}
Then there is
a partition of $I$ into $I_a$ and $I_b$ such that
$\sum_{i\in I_a}a_i>a-\Delta$ and $\sum_{i\in I_b}b_i\geq b$.
\end{lem}

\begin{proof}
Define a total order $\preceq$ on $I$ in a way that
$i\preceq j$ implies $\frac{a_{i}}{b_{i}}\leq\frac{a_{j}}{b_{j}}$ for all
$i,j\in I$. Let $\ell\in I$ be minimal in this order with $a\geq\sum_{i\succ \ell}a_{i}$.

Set $I_{a}:=\{i\in I:\ i\succ\ell\}$ and set $I_{b}:= I\setminus I_{a}$. It is
clear that $\sum_{i\in I_a}a_{i}>a-\Delta$, by the definition of $\ell$ and as
$a_{\ell}\leq \Delta$. So, all we have to show is that $\sum_{i\in I_b}b_{i}\geq b$.

Indeed, suppose otherwise. Then by (\ref{sum}), and by the definition of
$\ell$, we have that
\begin{align*}
\frac{\sum_{i\in I_b}b_{i}}{\sum_{i\in I}b_{i}} & <\ \frac{b}{\sum_{i\in I}b_{i}}\\
&\leq\ \frac{a-\sum_{i\in I_a}a_{i}}{\sum_{i\in I}a_{i}}+\frac{b}{\sum_{i\in I}b_{i}}\\
& \leq \
1-\frac{\sum_{i\in I_a}a_{i}}{\sum_{i\in I}a_{i}} \\ & =\ \frac{\sum_{i\in
I_{b}}a_{i}}{\sum_{i\in I}a_{i}}.
\end{align*}

Multiply the two sides of this inequality with $\sum_{i\in I}a_{i}\cdot
\sum_{i\in I}b_{i}$, subtract the term $\sum_{i\in I_b}a_{i}\cdot \sum_{i\in
I_b}b_{i}$, and divide by $\sum_{i\in I_a}b_i\sum_{i\in I_b}b_i$ to obtain
\[
\frac{a_\ell}{b_\ell}\leq \frac{\sum_{i\in I_a}a_{i}}{\sum_{i\in
I_a}b_{i}}<\frac{\sum_{i\in I_b}a_{i}}{\sum_{i\in I_b}b_{i}}\leq
\frac{a_{\ell}}{b_{\ell}},
\]
(where the first and last inequality follow from the definition of $\preceq$).
This yields the desired contradiction.
\end{proof}

We shall now apply Lemma~\ref{lem:numbers} to partition our matching $M=\{e_i\}_{i\leq |M|}$. We do this separately for the two cases from Section~\ref{prep}.

In Case~1, we set $$a:=|V(F_A)|+\frac{\pi k}{20},\  \ \ b:=|V(F_B)|+\frac{\pi
k}{20},\ \ \text{ and } \ \ \Delta:=2s.$$ For $i\leq |M|$, set $a_i:=\weightdeg_{e_i}(A) \leq \Delta$, and $b_i:=\weightdeg_{e_i}(B)\leq \Delta$. Now,~\eqref{eq:degCaseA}  implies that
\[
\frac{a}{\sum_{i=1}^{|M|}
a_i
 }+\frac{b}{\sum_{i=1}^{|M|}b_i}\leq \frac{|V(F_A)|+|V(F_B)|+\frac{\pi k}{10}}{(1+\frac{\pi}{10})k}\leq 1.
\]

Hence,  Lemma~\ref{lem:numbers} yields a partition of  $M$ into $M_A$ and $M_B$ such that

\begin{equation}\label{space1}
\weightdeg_{M_A}(A)>|V(F_A)|+\frac{\pi k}{40}\text{ \ and\    }\weightdeg_{M_B}(B)>|V(F_B)|+\frac{\pi k}{40}.
\end{equation}

\bigskip

In Case~2, set $$a:=|V(F_f)|+\frac{\pi k}{20},\ \ \ b:=|V(\bar F_B)|+\frac{\pi k}{40},\ \ \ \text{ and }\ \ \ \Delta:=2s.$$
For $i=1,\dots ,|M|$, again set $a_i:=\weightdeg_{e_i}(A) $, and $b_i:=\weightdeg_{e_i}(B)$. Set $L':=L\setminus M$. For $i=|M|+1,\dots ,|M|+|L'|$, set $a_i:=0$, and set $b_i:=\weightdeg_{C_i}(B)$, where $C_i$ is the $i$th cluster in $ L'$. 

Observe that by~\eqref{barA}, 
\begin{equation}\label{eq:Vb}
|V(\bar F_B)|\leq \frac {k-|V(F_f)|}{2}.
\end{equation}

  Now, let us check that the conditions of Lemma~\ref{lem:numbers} hold. Clearly, $a_i,b_i\leq \Delta$ for all $i\leq |M|+|L'|$.
  
Moreover, Condition~\eqref{sum}  holds since~\eqref{eq:degCaseB} and~\eqref{eq:Vb} imply that
 \begin{align*}
 \frac{a}{\sum_{i=1}^{|M|+|L'|}a_i}+\frac{b}{\sum_{i=1}^{|M|+|L'|}b_i} & \leq\ 
 \frac{|V(F_f)|+\frac{\pi k}{20}}{(1+\frac{\pi }{10})k}+\frac{|V(\bar
 F_B)|+\frac{\pi k}{40}}{(1+\frac{\pi }{10})\frac k2}\\ & =\  \frac{|V(F_f)|+2|V(\bar F_B)|+\frac{\pi k}{10}}{(1+\frac{\pi }{10})k}\\ & \leq \ 1.
 \end{align*}
 We thus obtain a partition of $M$  into $M_f$ and $\bar M_B$
 such that

\begin{equation}\label{space2}
\weightdeg_{M_f}(A)> |V(F_f)|+\frac{\pi k}{40}\quad 
\mbox{and}\quad\weightdeg_{\bar M_B\cup L'}(B)\geq |V(\bar F_B)|+\frac{\pi k}{40}.
\end{equation}

We partition  $\bar{\mathcal  T}_B$ into $\mathcal T^M_B\cup \mathcal T_B^L$ such that $\mathcal T^M_B$ will be
embedded using the edges of $\bar M_B$ and  $\mathcal T_B^L$ 
will be embedded using the clusters in $L'$. This
partition is necessary: we have to embed as much of
$\bar{\mathcal  T}_B$ as possible in the edges of $\bar M_B$, before we start
using the high average degree of clusters in $L'$, as the latter may alter the
possibility of using edges from $\bar M_B$.

 Let $\mathcal T^M_B\subseteq
\bar{\mathcal T}_B$ be maximal with

\begin{equation}\label{eq:degM_B2}
\weightdeg_{\bar M_B}(B)\geq |\bigcup_{T\in \mathcal T^M_B}V(T)|+\frac{\pi
k}{40N}|\bar M_B|.
\end{equation}

Set $\mathcal T_B^L:=\bar{\mathcal  T}_B\setminus \mathcal T_B^M$. Let
$F_B^M:=\bigcup_{T\in \mathcal T^M_B}T$ and let $F_B^L:=\bar F_B- V(F_B^M)$.

Observe that if $\mathcal T_B^M\neq \bar{\mathcal  T}_B$, then the maximality
of $\mathcal T_B^M$ ensures that $$\weightdeg_{\bar M_B}(B)<
|V(F^M_B)|+\frac{\pi k}{40N}|\bar M_B|+\beta k.$$ Hence,
by~\eqref{space2}, either $\mathcal T_B^L=\emptyset$, or
\begin{equation}\label{eq:deg_L'}
\weightdeg_{L'}(B)\geq |V(F_B^L)|+\frac{\pi k }{80N}|L'|.
\end{equation}

\subsection{Embedding lemmas for trees}\label{sec-embeddingLemma}

In this section, we shall prove some preparatory lemmas on embedding trees in regular pairs of $\bar H$.
As mentioned in the overview, it is important to keep the edges of
the matching in $\bar H$ \emph{balanced} as long as the edge is not
\emph{saturated}, i.\,e., as long as we did not embed in the regular pair the
expected number of vertices of the tree. This is captured below by property $(\star)$, where $U$ stands
for vertices already used in previous steps of the embedding
process, and $N$ stands for the neighbourhood of the image of the corresponding
seed mapped in cluster~$A$ or~$B$. So property $(\star)$ can be read as \emph{If the edge is not balanced, then it is
saturated}.

Let $C,D\in V(\bar H)$, and let $U,N\subseteq C\cup D$.
We say that $U$ has property~$(\star)$ in  $CD$ for
$N$ if it satisfies the following.
\begin{itemize}
  \item [$(\star)$] If $||C\cap U|-|D\cap U||>\beta k+\varepsilon s$, then \\
  $\min\{|N\cap C|,|N\cap D|\}\leq \min\{|C\cap U|,|D\cap U|\}+2\varepsilon s+\beta
  k$.
\end{itemize}

Now our first embedding lemma states that property $(\star)$ can be kept
throughout the embedding process. 

\begin{lemma}\label{lem:embeddingsmalltrees}
Let $T$ be a tree with root $r$ and of order at most $\beta k$. Let
$CD\in E(\bar H)$. 
Suppose that $U,N\subseteq C\cup D$ are such that
\begin{equation}\label{ncudu}
\min\{|N\cap C\setminus U|,|D\setminus U|\}> \frac 2p (\varepsilon s+\beta k).
\end{equation}
Then there is an embedding $\varphi$ of $T$ in $(C\cup
D)\setminus U$ such that $\varphi(r)\in N\setminus U$ and such that the
following holds.
\begin{itemize}
\item[$(\star\star)$] If $ U$ has property~$(\star)$ in  $CD$ for $N$, \\then  also $U_\varphi:= U\cup \varphi (V(T))$  has property~$(\star)$ in  $CD$ for $N$.
\end{itemize}
\end{lemma}

\begin{proof}
Write
$V(T)=r\cup L_1\cup L_2\cup\ldots$, where $L_\ell$ is the $\ell$th level of $T$ (i.\,e.\ the set of vertices at distance $\ell$ to $r$).

First, suppose that
$|N\cap D\setminus U|\leq \varepsilon s$.
In this case, choose $\varphi (r)\in N\cap C\setminus  U$ typical
to $D\setminus  U$. This is possible because by~\eqref{ncudu}, $|N\cap C\setminus U|>
\frac 2p(\varepsilon s+\beta k)>\varepsilon s$ and by~\eqref{eq:for regular pair},  at
most $\varepsilon s$ vertices of $C$ are not typical to the significant subset $D\setminus U$ of $D$.

Embed the rest of $V(T)$ levelwise.  For $\varphi (L_\ell)$, the image of the
$\ell$th level $L_\ell$, we choose unused vertices of $D\setminus  U$ that are
typical to $C\setminus  U$ if $\ell$ is odd, and unused vertices of
$C\setminus  U$ that are typical to $D\setminus  U$ if $\ell$ is even. Because
$C\setminus U$ and $D\setminus U$ are significant sets, any vertex that is
typical to $C\setminus U$, or to $D\setminus U$, has at least
$(p-\varepsilon)|C\setminus U|\ge \varepsilon s+\beta k$, resp.\ $(p-\varepsilon)|D\setminus U|\ge \varepsilon s+\beta k$, neighbours in $C\setminus U$, resp.\ in
$D\setminus U$ (here we used~\eqref{ncudu}). Among these
neighbours there are then at least $\beta k \geq V(T)$ vertices that are typical.  

Now, suppose that $|N\cap D\setminus  U|> \varepsilon s$. In this case, we may
alternatively wish to embed $r$ in $N\cap D$. We do so in either of the
following cases
\begin{enumerate}
\item $|\bigcup_{\ell \in\mathbb N}L_{2\ell-1}|>|\bigcup_{\ell\in\mathbb N}L_{2\ell}|$ and $|C\setminus  U|\geq |D\setminus  U|$, or
\item $|\bigcup_{\ell \in\mathbb N}L_{2\ell-1}|<|\bigcup_{\ell\in\mathbb N}L_{2\ell}|$ and $|C\setminus  U|\leq |D\setminus  U|$,
\end{enumerate}
and otherwise embed $r$ in $N\cap C$, as before. 
The purpose of embedding $r$ in~$D$ and not in $C$ is to keep the pair $(C,D)$
balanced, i.\,e., our choice of $r$ ensures that (if $|N\cap D\setminus  U|> \varepsilon s$)

\begin{equation}\label{diffatmostbetak}
||C\cap U_\varphi|-|D\cap U_\varphi||\leq \max\{||C\cap U|-|D\cap U||, \beta k\}
\end{equation}

Then, the rest of $T$ is embedded analogously as above
(possibly swapping the roles of $C$ and $D$). This completes the embedding of $T$.

\medskip

It remains to prove~$(\star\star)$. So assume that $U$ has property~$(\star)$ for $N$ in $CD$. Furthermore, assume that 
\begin{equation}\label{eq:big}
 ||C\cap U_\varphi|-|D\cap U_\varphi||>\beta k+\varepsilon s.
\end{equation}
 Now, if
$||C\cap  U|-|D\cap  U||>\beta k+\varepsilon s$, then property~$(\star)$
for $U_\varphi$ follows from property~$(\star)$ for $ U$.
Suppose otherwise, that is
\begin{equation}\label{j-2}
||C\cap  U|-|D\cap  U||\leq\beta k+\varepsilon s.
\end{equation}
By~\eqref{diffatmostbetak}, inequality~\eqref{eq:big} only holds if we could not choose
where to embed the root of $T$, in $N\cap C$ or in $N\cap D$. 
Hence,
$$|N\cap D\setminus U|\leq \varepsilon s.$$ Using~\eqref{j-2}, this gives
\begin{align*}
\min\{|N\cap C|,|N\cap D|\}
\leq \ &\max\{|C\cap  U|,|D\cap  U|\}+\min_{Y=C,D}\{|N\cap Y\setminus  U|\}\\
\leq \ &\max\{|C\cap  U|,|D\cap  U|\}+\varepsilon s\\
\leq \ &\min\{|C\cap  U|,|D\cap  U|\}+2\varepsilon s+\beta k\\
\leq \ &\min\{|C\cap U_\varphi|,|D\cap U_\varphi|\}+
2\varepsilon s+\beta k,
\end{align*}
as desired.
\end{proof}

We need some definitions.
Let $C,D,X\in V(\bar H)$, 
We say that $U\subseteq V(G)$ has property~$(\diamond)$ in  $(C,D)$ with respect to
$X$ if it satisfies the following.
\begin{itemize}
  \item [$(\diamond)$] If $||C\cap U|-|D\cap U||>\beta k+\varepsilon s$, then \\
  $\min\{\weightdeg_C(X),\weightdeg_D(X)\}\leq \min\{|C\cap U|,|D\cap
  U|\}+4\varepsilon s+\beta k$.
\end{itemize}

Let $X'\subseteq X$, let $v\in X$, let $\mathcal Z\subseteq V(\bar H)$.
 An embedding $\varphi$ of a rooted tree  $(T,r)$ is a {\em
$(v,X',U)$-embedding} in $\mathcal Z$,  if $\varphi(V(T)\setminus \{r\})\subseteq \bigcup \mathcal Z\setminus U$, if $\varphi(r)=v$, and if each
vertex at odd distance to the root $r$ is mapped to a vertex that is typical to
$X'$. 
 A vertex is {\em $\mathcal Z$-typical}, if it is typical to each cluster
from~$\mathcal Z$. 
For each cluster $C\neq X$, let $C_{X'}$ be the set of all vertices of $C$ that are not typical to $X'$, and let $S_{X'}:=\bigcup_{C\in V(\bar H), C\neq X}C_{X'}$. Note that $C_{X'}=\emptyset$ if $d(C,X)=0$.

Finally, for
$m\in \mathbb N$, the set  $\mathcal Z$ is said to be {\em $(m,U)$-large for $X$}, if
\begin{equation*}
\weightdeg_{\mathcal Z}(X)>m+|U\cap \bigcup \mathcal Z|+\frac{\pi k}{100
N}|\mathcal Z|.
\end{equation*}

\begin{lemma}\label{lem:forestInMatching}
Let $T,r, X',X,v$ and $U$ be as above with
$|X'|\geq |X|/2$. \\ 
 {\bf A)}~Suppose $M_X$ is a matching in
$\bar H-X$ so that $V(M_X)$ is
$(|V(T)|,U)$-large for~$X$, so that  $v$ is $V(M_X)$-typical, and so that
$U\cup S_{X'}$ has property~$(\diamond)$ in $(C,D)$ with respect to $X$, for each $CD\in
M_X$. \\
Then, there is a $(v,X',U)$-embedding
 $\varphi$ of $T$ in  $V(M_X)$ such that $U\cup
\varphi(V(T))\cup S_{X'}$ has property~$(\diamond)$ with respect to $X$ for every $CD\in
M_X$.\\
{\bf B)}~Let $L_X, W_X \subseteq V(\bar H)$ be such that $L_X$ is
$(|V(T)|,U)$-large for $X$,  and $W_X$ is $(|V(T)|,U)$-large for each $Y\in
L_X$. If $v$ is $L_X$-typical, then there is a $(v,X',U)$-embedding
$\varphi$ of $T$ in $L_X\cup W_X$.
\end{lemma}
 
\begin{proof}
We map $r$ to $v$ and embed the trees from the forest $F:=T-\{r\}$
inductively. In each step $j\geq 1$, we embed a tree $T^j$ of the forest $F$. 
Denote by $V^{j}$ the set $\bigcup_{i\leq j}V(T^i)$ of
vertices we have embedded just after step $j$ and set $V^0=\emptyset$. 
Set
$U^{j}:=U\cup S_{X'}\cup\varphi (V^{j})$ for any $j\ge 0$. In particular, $U^{0}=U\cup S_{X'}$.

For Part~{A)}, we shall ensure the following two properties of $U$ during our
embedding. Firstly, 
 if $CD\in M_X$ satisfies $ | |C\cap U^0|-|D\cap U^0| |\leq 
\beta k+\varepsilon s$, then we require that for every $j\geq 1$
\begin{itemize}
  \item [(I)]$U^{j-1}{\rm ~has~property~(}\star{\rm )~for~}N(v).$
\end{itemize}
This property
holds for $j=1$, as the condition of property~$(\star)$ is void, and we shall check it for each later step.

 Secondly, for those edges with $ | |C\cap U^0|-|D\cap U^0| |> \beta k+ \varepsilon s$, observe that as
the sets $U^{j}$ are growing, property~$(\diamond)$ ensures that for all $j\geq 1$
\begin{itemize}
  \item [(II)]$\min_{Y\in \{C,D\}}\{\weightdeg_Y(X)\}\leq \min_{Y\in
  \{C,D\}}\{|Y\cap U^{j-1}|\}+4\varepsilon s+\beta k.$
\end{itemize}

So, assume now that we are in step $j\geq 1$, that is, $\varphi(x)$ has been
defined for all $x\in V^{j-1}$, and we are about to embed $T^j$.

\begin{claim}\label{obs:edge-existence}
There is an edge $CD$, with $CD\in M_X$ for Part { A)} and  with $C\in
L_X$, and $D\in W_X$ for Part { B)}, such that $$\min\{|(N(v)\cap C)\setminus
U^{j-1}|,|D\setminus U^{j-1}|\}\geq \frac 2p (\varepsilon s + \beta k).$$
\end{claim}
Before proving Claim~\ref{obs:edge-existence}, we shall show how we complete our embedding
of $T^j$ under the assumption that the claim holds for some edge $e:=CD$.

Set $N:=N(v)\cap e$ and let $r^j:=N(r)\cap V(T^j)$ be the root of $T^j$. 
Use Lemma~\ref{lem:embeddingsmalltrees} 
 to embed $T^{j}$ in $e\setminus U^{j-1}$, mapping $r^j$ to
 $N\setminus U^{j-1}$. Lemma~\ref{lem:embeddingsmalltrees} together with~(I)
 for $j$ ensures (I) for $j+1$. As our embedding avoids $S_{X'}$, all vertices in $\varphi(T^j)$ are typical to $X'$. 
This terminates step $j$. 

Say we terminate the embedding procedure after step $\ell$ (that is, $\ell$ is the number of components of $F$). Then $\varphi$ is a
$(v,X',U)$-embedding. So, for Part B), we are done. For Part A), however, we still have to prove that
$U\cup \varphi(V(T))\cup S_{X'}$ 
has property~$(\diamond)$ in $(C,D)$ with respect to $X$,  for each $CD\in M_X$.

To this end, assume that 
\begin{equation}\label{eq:endLemma11}
||C\cap U^\ell|-|D\cap U^\ell||>\beta k+\varepsilon s\;.
\end{equation}
If $||C\cap U^0|-|D\cap U^0||\leq\beta k+\varepsilon s$, then~(I) holds  by induction
for $\ell+1$ and thus
$U^{\ell}$  has property~$(\star)$
in $CD$ for $N(v)$.
Hence, because $v$ is typical to $C$ and $D$,
\begin{align*}
\min_{Y=C,D}\{\weightdeg_Y(X)\}&\leq
\min_{Y=C,D}\{\deg_Y(v)\}+\varepsilon s\\ 
&\overset{\eqref{eq:endLemma11},(\star)}{\leq} \min_{Y=C,D}\{|Y\cap U^{\ell}|\}+3\varepsilon s+\beta k\\
&\leq \min_{Y=C,D}\{|Y\cap (U^{\ell}\setminus S)|\}+4\varepsilon s+\beta k\\
&= \min_{Y=C,D}\{|Y\cap (U\cup V(T))|\}+4\varepsilon s+\beta k.
\end{align*}
On the other hand, 
  if $||C\cap U^0|-|D\cap U^0||>\beta k+\varepsilon s$,  
then~(II) ensures that 
 $U^\ell|=U\cup \varphi(V(T))\cup S_{X'}$ has property~$(\diamond)$
in each $CD\in M_X$ for Part~{ A)}. It only remains to
prove Claim~\ref{obs:edge-existence}.
\bigskip

{\bf Proof of Claim~\ref{obs:edge-existence}:}
First, suppose we are in Case {A)}. Let us start by showing that there is an
edge $e=CD\in M_X$ which satisfies
\begin{equation}\label{edge!}
\weightdeg_e(X)-|e\cap U^{j-1}|\geq \frac {8}p(\varepsilon s +\beta
k)+2\varepsilon s.
\end{equation}
Indeed, suppose there is no such edge. Then, as $V(M_X)$ is $(|V(T)|,U)$-large,
we have that
\begin{align*}
\frac{8}p(\varepsilon s+\beta k)|M_X|>& \sum_{e\in M_X}
(\weightdeg_{e}(X) -|e\cap U^{j-1}|-2\varepsilon s)\\
=\ &\weightdeg_{M_X}(X)-|U\cap \bigcup M_X|-|U^{j-1}\setminus U|-2\varepsilon
s|M_X|\\ \geq\  & \weightdeg_{M_X}(X)-|U\cap \bigcup
M_X|-|V(T)|-|S_{X'}\cap M_X|-2\varepsilon s|M_X|\\ \geq \ & \frac{\pi k}{100
N}|V(M_X)|-2\varepsilon s|M_X|\\ >\ &\frac{\pi k}{100N}|M_X|,
\end{align*}
which, as $\beta k\leq \frac {\varepsilon}{M_0}n\leq \varepsilon s$, implies that $16\varepsilon/p>\pi q/100$, a contradiction.

So, assume now that we have chosen an edge $e$ for which~\eqref{edge!} holds. 
Clearly, we can write $e=CD$ such that
\begin{align}\label{freeC!}
 \frac 4p(\varepsilon s+\beta k)&\overset{\eqref{edge!}}{\leq}
 \weightdeg_C(X)-\varepsilon s-|C\cap U^{j-1}|\\
 &\leq |N(v)\cap C\setminus  U^{j-1}|.\label{freeC!!}
\end{align}
We claim that 
\begin{equation}\label{D}
|D\setminus U^{j-1}|\geq \frac 2p(2\varepsilon s +\beta k),
\end{equation}
 which
together with~\eqref{freeC!!} implies Claim~\ref{obs:edge-existence} for Case {A)}.
Indeed, suppose for contradiction~\eqref{D} does not hold. Then~\eqref{freeC!}
implies that
\begin{align}\label{C-less-than-D!}
|C\cap U^{j-1}|& \leq \ s-\frac 4p(\varepsilon s+\beta k)-\varepsilon s \nonumber\\
& =\ |D\cap U^{j-1}|+|D\setminus U^{j-1}|-\frac 2p(2\varepsilon s+\beta k)-\frac
2p\beta k-\varepsilon s \nonumber\\ & <\ |D\cap
U^{j-1}|-\frac 2p\beta k -\varepsilon s.
\end{align}

We claim that
\begin{equation}\label{eq:byIorII}
\min\{\weightdeg_C(X),\weightdeg_D(X)\}\leq |C\cap U^{j-1}|+4\varepsilon s+\beta
k.
\end{equation}
Indeed, if $||C\cap U^0|-|D\cap U^0||\le \beta k+\varepsilon s$, then by~(I), $U^{j-1}$ has
property $(\star)$ for $N(v)\cap (C\cup D)$. As~\eqref{C-less-than-D!} implies that $||C\cap U^{j-1}|-|D\cap U^{j-1}|>\beta k +\varepsilon s$, we
obtain that 
\begin{align*}
\min\{\weightdeg_C(X),\weightdeg_D(X)\}
&\leq \min\{|N(v)\cap C|,|N(v)\cap D|\}+\varepsilon s\\
& \overset{(\star)}{\le} \min\{|C\cap U^{j-1}|,|D\cap
U^{j-1}|\}+3\varepsilon s +\beta k,
\end{align*}
 implying~\eqref{eq:byIorII}. On the other hand, if $||C\cap
U^0|-|D\cap U^0||> \beta k+\varepsilon s$, then~\eqref{eq:byIorII} follows
directly from~(II).

Thus, by~\eqref{edge!},
\begin{align*}
\frac {8}p(\varepsilon s+\beta k)+2\varepsilon s\leq &\ \weightdeg_e(X)-|C\cap
U^{j-1}|-|D\cap U^{j-1}|\\ \overset{\eqref{eq:byIorII}}{\leq} &\
\weightdeg_e(X)-\min_{Y\in
\{C,D\}}\{\weightdeg_Y(X)\}+4\varepsilon s+\beta k-|D\cap
U^{j-1}|\\ \leq &\ s+4\varepsilon s+\beta k-|D\cap U^{j-1}|\\ < &\ |D\setminus U^{j-1}|+4\varepsilon s+\beta k.
\end{align*}
So, $|D\setminus U^{j-1}|>(\frac 8p-2)(\varepsilon s+\beta k)$, a contradiction
to our assumption that~\eqref{D} does not hold. This proves~\eqref{D}.

\medskip 

Now, assume that we are in Case {B)}. 
First we show that if some $\mathcal Z\subseteq V(\bar H)$ is
$(|V(T)|,U)$-large for some $Y\in V(\bar H)$, then there is a $Z\in \mathcal
Z$ such that 
\[
\weightdeg_{Z}(Y) -|Z\cap  U^{j-1}|\geq \frac 2p(\varepsilon s+\beta
k)+\varepsilon s,
\] 
which implies that $Z\in N(Y)$.

Indeed, otherwise, by the definition of $(V(T),U)$-large and using the fact that 
$|V(T)|+|U\cap \bigcup \mathcal{Z}|\ge |U^{j-1}\cap \bigcup \mathcal Z|-\varepsilon s|\mathcal Z|$, we have that
\begin{align*}
\frac 2p(\varepsilon s+\beta k)|\mathcal Z|&>\ \sum_{Z\in
\mathcal Z}(\weightdeg_{Z}(Y)-|Z\cap U^{j-1}|-\varepsilon s)\\
&=\ \weightdeg_{\mathcal Z}(Y)-|U^{j-1}\cap \bigcup \mathcal
Z|-\varepsilon s|\mathcal Z|\\
&>\ (\frac {\pi k}{100 N}-2\varepsilon s)|\mathcal Z|\\
&\geq\ \frac {\pi k}{200 N}|\mathcal Z|,
\end{align*}
a contradiction. 

Applying this assertion with $\mathcal Z=L_X$ and $Y=X$, we obtain $C\in L_X$
such that \[|N(v)\cap C\setminus U^{j-1}|\geq \weightdeg_C(X)-|C\cap
U^{j-1}|-\varepsilon s\geq \frac 2p(\varepsilon s+\beta k).\]
Applying the assertion again with $\mathcal Z=W_X$ and $Y=C$, we obtain $D\in
W_X\cap N(C)$ such that 
\[|D\setminus
U^{j-1}|\geq \weightdeg_D(C)-|D\cap U^{j-1}|\geq \frac 2p(\varepsilon s+\beta
k),\] as desired for Claim~\ref{obs:edge-existence}.


\end{proof}

\subsection{The embedding in Case~1}\label{sec:1}

In this subsection, we shall complete the proof of Theorem~\ref{thm:approx}
under the assumption that Case~1 of Section~\ref{prep} holds. So, we assume
that there are an edge $AB\in E(\bar H)$  and a matching $M=M_A\cup M_B$ in 
$\bar H-\{A,B\}$ as in Section~\ref{num-part}. These, together with
the  sets $SD=SD^A\cup SD^B$, $F_A$ and $F_B$ from Section~\ref{part},
satisfy~\eqref{space1}.

\smallskip

Our embedding $\varphi$ will be defined in $|SD|$ steps.
In each step $i\geq 1$, we choose a suitable vertex $r_i\in SD$ and embed it
together with all trees from  $$\mathcal T_i:=\{T\in \mathcal T_A\cup \mathcal T_B:\ sd(T)=r_i\}.$$

Set $V_0:=\emptyset$ and for $i\geq1$, let $$V_{i}:=V_{i-1}\cup\{r_i\}\cup
\bigcup_{T\in \mathcal T_i}V(T).$$

We start with the root $r_1:=R$ of $T^*$, and in each step $i>1$, we shall
choose a vertex $r_i\in SD\setminus V_{i-1}$ that  is adjacent to $V_{i-1}$.  The seed
$r_i$ will be embedded in a vertex $v_i\in A\cup B$, while $\mathcal T_i$ will be mapped to
edges from $M$ (or more precisely, to the corresponding subgraph of $G_p$).
 Set $U_0:=\emptyset$, and once $\varphi$ is defined on $V_i$, set $U_i:=\varphi (V_i)$. 

For each $i\geq 0$, the following conditions will hold.

\begin{enumerate}[(i)]
\item \label{FewInA,B} $|(A\cup B)\cap U_i|\leq i$,

\item \label{typicalInC} if $x\in V_i\cap N(SD^A)$, resp.\ $x\in V_i\cap N(SD^B)$, then $\varphi(x)$ has at least
  $\frac p4s$ neighbours 
 in $A$, resp.\ in $B$,

\item \label{well-embeddedA!NEW!} for $CD\in M_A$, the set $U_i\cup S_A$ has
property~$(\diamond)$ in $CD$ with respect to $A$. 

\item \label{well-embeddedB!NEW!} for $CD\in M_B$, the set $U_i\cup S_B$ has
property~$(\diamond)$ in $CD$ with respect to $B$. 
\end{enumerate}
Observe that properties~(i)--(iv) trivially hold for $i=0$.

\medskip

So, suppose now that we are in some step $i\geq 1$ of our embedding process. Choose~$r_i\in SD$ as detailed above.
Let us assume that $r_i\in SD^A$, the case when $r_i\in SD^B$ is analogous.

We embed $r_i$ in a vertex  $
v_i=\varphi(r_i)\in A$ that is typical to~$B$ and typical to all but at most
$2\sqrt{\varepsilon}|M_A|$ clusters of $M_A$. Properties~\eqref{FewInA,B} and~(\ref{typicalInC}) for $i-1$ ensure
that if $x$ is the predecessor of $r_i$ in $T^\ast$, then~$\varphi(x)$ has at
least $\frac {ps}{4}-i$ neighbours in $A\setminus U_{i-1}$.
By~\eqref{eq:for regular pair} and~\eqref{eq:typicality2}, at most
$2\sqrt{\varepsilon} s$ of these vertices do not have the required properties.
Hence, there are at least $(\frac p4-2\sqrt{\varepsilon})s-i\geq 1$ suitable vertices we may choose $v_i$ from.

Let $M_A^i\subseteq M_A$ be a maximal
submatching such that $v_i$ is typical to each of the end-clusters of
each edge of $M_A^i$, i.\,e., $v_i$ is $V(M_A^i)$-typical. Then
by~\eqref{eq:typicality2} and~\eqref{space1} we obtain
\begin{align}\label{eq:M_A^i}
\weightdeg_{M_A^i}(A)&\geq \ \weightdeg_{M_A}(A)-4\sqrt{\varepsilon}  |M_A|s\notag \\
&>\ |V(F_A)|+\frac{\pi k}{40}-4\sqrt{\varepsilon}  Ns\notag\\
&>\ |V(F_A)|+\frac{\pi k}{80}\notag\\
&>\ |\bigcup_{T\in  \mathcal T_i}V(T)|+|U_{i-1}\cap
\bigcup_{C\in V(M_A)}C |+\frac{\pi k}{80N}|V(M_A^i)|.
\end{align}
 
Let $T$ be the tree induced by $r_i$ and the trees from $\mathcal T_i$, and let
$r:=r_i$  be the root of $T$. Each component of $T-r$ has order at most $\beta
k$. Inequality~\eqref{eq:M_A^i} implies that $V(M_A^i)$ is $(|V(T)|,U_{i-1})$-large
for $A$. 
Observe that $U_{i-1}\cup S_A$ has property $(\diamond)$ in $(C,D)$ with
respect to $A$ for each $CD\in M_A^i$ by~\eqref{well-embeddedA!NEW!}.

Now we use Lemma~\ref{lem:forestInMatching} Part~{A)} with $T$ and
setting $M_X:=M_A^i$, $U:=U_{i-1}$,  $v:=v_i$, and $X=X'=A$. This provides
with a $(v_i,A,U_{i-1})$-embedding of $T$ in $V(M_A^i)$. Thus every
vertex of $T-r$ at odd distance from $r$ is mapped to a vertex that
is typical to~$A$, i.\,e., that has at least $(p-\varepsilon)|A|\geq
\frac p4s$ neighbours in~$A$. By (II) and (VI) of Section~\ref{part} this
implies that \eqref{typicalInC} holds for all vertices in $V(T-r)\cap
N(SD)$. For $r$ property~\eqref{typicalInC} is satisfied as
$v_i$ is typical to~$B$ and thus has at least $(p-\varepsilon)|B|\geq
\frac p4s$ neighbours in~$B$.  It is easy to see that~\eqref{FewInA,B} holds for $i$,
as it holds for $i-1$, and by our choice of $\varphi(V_i\setminus V_{i-1})$.
Property~\eqref{well-embeddedB!NEW!} trivially holds as no vertices were mapped to~$M_B$.
Lemma~\ref{lem:forestInMatching} Part~{A)} ensures property~$(\diamond)$ for
all edges $CD\in M_A^i$. Because we did not embed anything in the edges of 
$M_A\setminus M^i_A$,~\eqref{well-embeddedA!NEW!} for $i-1$
implies~\eqref{well-embeddedA!NEW!} for $i$, for all $CD\in M_A$.

This completes the embedding of the tree $T^\ast$ in $G_p\subseteq G$ in
Case~1.

\subsection{The embedding in Case~2}\label{sec:2}

We shall now complete the proof of Theorem~\ref{thm:approx}
under the assumption that Case~2 of Section~\ref{prep} holds. That is, there
are an edge $AB\in E(\bar H)$ and a  matching $M= M_f\cup\bar M_B$ in
$\bar H-\{A,B\}$ together with sets $\SD=\SD^A\cup SD^B$, $F_f$,
$\bar{F}_A$, $F_B^M$ and $F_B^L$ from Sections~\ref{part}
and~\ref{switch} satisfying ~\eqref{space2},~\eqref{eq:degM_B2}
and~\eqref{eq:deg_L'} from Section~\ref{num-part}.

\smallskip

Our embedding will be defined in three phases. In the first phase, we shall
embed all vertices from $\SD$ in $A\cup B$, embed $F_f$ in  edges of
$M_f$, and embed~$F_B^M$ in edges of $\bar M_B$. In the second phase, we shall
embed $F_B^L$ in edges incident with $L'\cap N(B)$, and in the third phase, we
shall embed  $\bar F_A$ in the remaining space inside edges from $M$.

Denote by $A'$ the set of vertices in $A$ that are typical to all but at
most $2\sqrt{\varepsilon}|M|$ clusters of $V(M)$, and denote by $B'$ the set of
vertices in $B$ that are typical to all but at most $\sqrt{\varepsilon}|L'|$
clusters of $L'$.

The first phase is done analogously as in Case~1, while considering $A'$
and $B'$ instead of $A$ and $B$. In
each step, Lemma~\ref{lem:forestInMatching} Part~{A)} is used in the following setting.

The tree $T$ is 
the tree induced by $r_i$ and the trees from $$\mathcal T_i:=\{T\in \mathcal
T_f\cup \mathcal T_B^M:\; sd(T)=r_i\}.$$ Its root is $r:=r_i$. We set either
$(X',X)=(A',A)$ or $(X',X)=(B',B)$, and let $v=\varphi(r_i)$.  
The matching $M_X$ is a maximal submatching either of $M_f$ or of $\bar
M_B$, so that $\varphi(r_i)$ is $V(M_X)$-typical. Finally, the set $U$ is the set of
the vertices used before step~$i$.

\medskip

For the second phase, assume that $V(F_B^L)\neq\emptyset$ (otherwise we shall skip the second phase).
We define the second phase of our embedding process in $| {SD}^B|$ steps.

In each step $i\geq 1$, we embed the trees $\mathcal T^i:=\{T\in \mathcal
T_B^L\::\: sd(T)=r_i\}$ in edges incident with $L'$. (Recall that $L'=L\setminus M$.) Suppose that we are at step
$i$ of this procedure, i.\,e.~that we have already embedded the trees
from $\mathcal T^1,\dots,\mathcal T^{i-1}$. Denote by $
U_{i-1}$ the set of vertices used so far for the embedding. 
Let $L'_i$ be the set of those clusters of $L'$ to which $\varphi(r_i)$ is typical.
As $\varphi(r_i)\in B'$,~\eqref{eq:typicality2} and~\eqref{eq:deg_L'} imply that 
$$\weightdeg_{L'_i}(B)\geq |\bigcup_{T\in \mathcal
T^i}V(T)|+|U_{i-1}\cap L'_i|+\frac {\pi k}{100 N}|L'_i|.$$
Furthermore, by~\eqref{Ldeg}, for all $Y\in L'_i$ we have that
$$ \weightdeg(Y)\geq |\bigcup_{T\in \mathcal
T^i}V(T)|+|U_{i-1}|+\frac
{\pi k}{100}.$$

Use Lemma~\ref{lem:forestInMatching} Part~{B)} to embed $\mathcal T_i$, letting the tree be the tree induced by $r_i$ and the trees from
$\mathcal T^i$, its root be $r_i$, and setting $X:=B$,  $X':=B'$,
$v:=\varphi(r_i)$, $L_X:=L'_i$, $W_X:=N(L'_i)$, and $U:=U_{i-1}$.

\medskip

The third phase of our embedding process takes place in $|\SD^A|$
steps, where in each step $i\geq 1$, we embed the trees from
$\mathcal T^i:=\{T\in \bar{\mathcal T}_A\::\: sd(T)=r_i\}$. 
Suppose that we are at step
$i$ of this procedure, i.\,e.~that we have already embedded the trees
from $\mathcal T^1,\dots,\mathcal T^{i-1}$. Denote by $
\bar U_{i-1}$ the set of vertices used so far for the embedding.
Let $M_i$ be the maximal submatching of $M$ such that $\varphi(r_i)$ is typical
to all cluster of $V(M_i)$. As $\varphi(r_i)\in A'$, we have by~\eqref{eq:typicality2} and~\eqref{eq:degCaseA} that
$$\weightdeg_{M_i}(A)\geq |V(\bigcup \mathcal T^i)|+|\bar U_i|+\frac {\pi
k}{100}.$$ 

Observe that, as each edge $CD\in M$ meets $N(A)$ in at most one end-cluster, the
set $U_i$ trivially has property~$(\diamond)$ in $CD$ with respect to $A$. We
use Lemma~\ref{lem:forestInMatching} Part~{A)} to embed $\mathcal T_i$, letting $T$ be the tree induced by $r:=r_i$ together with the trees
from $\mathcal T^i$, and setting $X:=A$,  $X':=A'$,
$v:=\varphi(r_i)$, $M_X:=M_i$, and $U:=\bar U_{i-1}$.

This terminated our embedding of $T^*$, and thus the proof of
Theorem~\ref{thm:approx}.

\section{Extensions and applications}\label{ext}

In this last section, we explore applications and generalisations of Theorem~\ref{thm:approx}. In Section~\ref{sec:ramsey} we show how our theorem implies an asymptotic upper bound on the Ramsey number of trees. We extend Theorem~\ref{thm:approx} so that it allows for embedding subgraphs other than trees in Section~\ref{sec:bip}.

\subsection{A bound on the Ramsey number of trees}\label{sec:ramsey}
Recall that $r(\mathcal H,\mathcal H')$ denotes the Ramsey number for the classes $\mathcal H$ and $\mathcal H'$ of graphs, and that $\mathcal T_{\ell}$ denotes the class of trees of order $\ell$.

Based on ideas from~\cite{discrepency} and using Theorem~\ref{thm:approx}, we
prove Proposition~\ref{rams:number}, which stated that $r(\mathcal T_{k+1},\mathcal T_{m+1})\leq k+m+o(k+m)$.
The sharp bound $k+m$ has been conjectured in~\cite{discrepency}.

\begin{proof}[Proof of Proposition~\ref{rams:number}]

Given $0<\varepsilon<1/4$, we apply Theorem~\ref{thm:approx} to
$\eta=q=\varepsilon /4$ to obtain an $n_0\in \mathbb N$. Now, let $n\geq n_0$, and let $G$ be a graph on $n'=(1+2\varepsilon)n+1$ vertices. Let $k$ and $m$ be such that $k+m=n$.

Clearly, either at least 
half of the vertices of $G$ have degree at least $k+\varepsilon n$, or in
the complement $\bar{G}$ of $G$, at least half of the vertices have degree at
least $m+\varepsilon n$.

First, suppose that the former of these assertions is true. Then it is easy to calculate that
\[
k+\varepsilon n\geq (1+\eta)(k+qn').
\]
Thus, we may apply Theorem~\ref{thm:approx}, which yields that each tree in $\mathcal T_{k+qn'+1}$ is a subgraph of $G$. Hence, also each tree in $\mathcal T_{k+1}$ is a subgraph of $G$.

Now, assume that the second assertion from above holds.
We have thus shown that for every $\varepsilon >0$ there is an
$n_0$ so that for all $k,m$ with $k+m\geq n_0$, we have that
$r(\mathcal T_{k+1},\mathcal T_{m+1})\leq (1+2\varepsilon)(k+m)+1$. This proves Proposition~\ref{rams:number}.
\end{proof}

\subsection{Graphs with few cycles}\label{sec:bip}

The question we pursue in this subsection is whether the condition of
Theorem~\ref{thm:approx} allows for embedding other graphs on $k+1$ vertices, apart from trees. For instance, may we add an edge to our tree $T^\ast$ and still embed it in $G$? In Theorem~\ref{thm:bip} we show that we may indeed add constantly many edges, as long as our graph stays bipartite.

Observe that the argument for the bound on
Ramsey number from Subsection~\ref{sec:ramsey} would apply here as well. We thus get an upper bound of $k+m+o(k+m)$ for the Ramsey numbers of graphs $Q_k$, $Q_m$ as in Theorem~\ref{thm:bip}, although the sharp bound does not hold (cf.~the example given in the introduction).

Our proof of Theorem~\ref{thm:bip} follows closely the lines of the proof of Theorem~\ref{thm:approx}. We embed a spanning tree $T^\ast$ of $Q$, and choosing $\varphi$ carefully, we ensure the adjacencies for the edges from $E(Q)\setminus E(T^\ast)$.

\begin{proof}[Proof of Theorem~\ref{thm:bip}]
Set $\pi:=\min\{\eta,q\}$ and set
$$\quad \varepsilon':=\frac {\varepsilon^{c+1}}{(c+3)^2}
, \quad\text{ and }\quad m_0:= \frac{500}{\pi^2 q},$$
where $\varepsilon$ is the constant from the proof of Theorem~\ref{thm:approx}.
As in the proof of Theorem~\ref{thm:approx},  the regularity lemma applied to
$\varepsilon'$, and $m_0$, yields natural numbers~$N_0$ and $M_0'$. Set
$M_0:=\max\{M'_0,c\}$, define $\beta$ and $p$ accordingly, and set
$$n_0:=\max\left\{N_0, \frac {9M_0}{\beta}\left(\frac {8}{p}\right)^{c+1}
\right\}.$$

Now, let $G$ be a graph on $n\geq n_0$ vertices which satisfies the condition of Theorem~\ref{thm:bip}, let $k\geq qn$, and let $Q$ be a connected bipartite graph of order $k+1$ with at most~$k+c$ edges, with a spanning tree $T^\ast$.  Fix a root $R$ in $T^\ast$.
Denote by $Q'$ the subgraph of $Q$ induced by the edges in $E(Q)\setminus
E(T^\ast)$ and let $P$ be the set of predecessors of $V(Q')$ in the tree order
of $T^\ast$.

We decompose~$T^\ast$ as in Section~\ref{part}, with the difference that we now
add the vertices from $V(Q')\cup P$ to the sets $A'$ and $B'$ (from the
definition of $SD$), depending on the parity of their distance in $T^\ast$ to
$R$. In this way, and since $Q$ is bipartite, we obtain, after the switching,
two independent sets~$\SD^ A$ and $SD^B$ so that $$|\SD^A|+ |SD^B|\leq \frac 8\beta+8c<\frac 9\beta,$$ which is constant in $n$.

The definition of our the embedding $\varphi$ is similar as in the proof of
Theorem~\ref{thm:approx}, except for some extra precautions we take for
vertices from $V(Q')\cup P$. At step~$i$, for each vertex $r\in \SD^A$, define
$$N_r^i:=\bigcap_{\ell=1}^{j}N(\varphi(x_\ell))\cap A,$$ where $x_1,\dots x_j$
are the already embedded neighbours of $r$ in $\SD^B$. If none of the neighbours of $r$ in $\SD^A$ has been embedded before step~$i$, then set~$N^i_r:= A$. Analogously define $N_r^i$ for $r\in SD^B$.

In each step $i$ of our embedding process, we shall ensure the following.

\begin{equation} \label{commonNeighbour}
\text{If }r\in V(Q')\text{ is not yet embedded, then
}|N_r^i|\geq \left(\frac p4\right)^js,
\end{equation}

where $j=j(r,i)$ is the number of neighbours of $r$ in $\SD^ A$ resp.~$SD^B$
that have already been embedded before step $i$.

Observe that
in step $i=0$, either $N_r^0=A$ or $N_r^0=B$, and
therefore~\eqref{commonNeighbour} is satisfied.

Suppose that at step $i\geq 1$ of our embedding process we are about to embed
a vertex $r=r_i\in V(Q')\cup P$. Assume that $r\in \SD^A$ (the
case when $r\in SD^B$ is analogous). Denote by $x_1,\dots, x_\ell$ the
neighbours of $r$ in $V(Q')$ that have not been embedded yet.

Now, embed $r$ in a vertex $v$ from $N_r^{i-1}$ that satisfies the three
following conditions of typicality:
\begin{itemize}
  \item $v$ is typical to all but at most $2\sqrt{\varepsilon}|M|$
  clusters of $V(M)$, resp.\ all but at most $\sqrt{\varepsilon}|L'|$ clusters
  of $L'$,
  \item $v$ is typical to all but at most $2\sqrt{\varepsilon}|M'|$
  clusters of the matching $M'$, where $M'$ stands either for $M_A$, $M_B$,
  $M_f$, or $\bar M_B$, depending on the case, and
  \item $v$ is typical to each $N_{x_j}^{i-1}$, for $1\leq j\leq
  \ell$.
\end{itemize}
This is possible, since our embedding
scheme and the condition on the number of edges of $Q$ ensure that $r$ has at
most $c+1$ neighbours in $Q$ that are already embedded. Thus, by~\eqref{commonNeighbour} for~$i-1$ and for~$r$, by~\eqref{eq:for regular
pair} and~\eqref{eq:typicality2}, and by choice of~$n_0$, there are
at least

 $$\left(\left(\frac
 p4\right)^{c+1}-(c+1)\varepsilon'-2\sqrt{\varepsilon'}\right)s-|\SD|+1\geq
 \frac 12\left(\frac p4\right)^{c+1}s-\frac 9\beta +1\geq 1$$unused typical vertices we can choose $\varphi(r)$ from.

Finally, observe that since we chose $\varphi(r)$ typical
to each~$N_{x_j}^{i-1}$, we have ensured
property~\eqref{commonNeighbour} for $i$ and for every $r'\in V(Q')$ that is
not yet embedded. This completes the proof of Theorem~\ref{thm:bip}.
\end{proof}

\section*{Acknowledgment}
The authors would like to thank Martin Loebl for helpful discussions and
Mikl\'os Simonovits for his valuable comments. We are also grateful to one of the referees for his/her
careful reading and numerous relevant remarks that helped to improve the
presentation of the paper considerably.


\thebibliography{llll}


\bibitem{aks}{M.~Ajtai, J.~Koml\'os and E.~Szemer\'edi,
On a conjecture of Loebl, In {\em Proc.\ of the 7th International Conference
on Graph Theory, Combinatorics, and Algorithms},  pages 1135-1146, Wiley, New
York, 1995.
}

\bibitem{blw}{C.~Bazgan, H.~Li and M.~Wo\'zniak,
On the Loebl-Koml\'os-S\'os conjecture, {\em J.~Graph Theory},
34:269-276, 2000.
}


\bibitem{bradob}{
S.~Brandt and E.~Dobson,
The Erd\H os--S\'os conjecture for graphs of girth $5$,
{\em Discr. Math.}, 150: 411-414,1996.
}

\bibitem{cooley}{
O.~Cooley, Proof of the Loebl-Koml\'os-S\'os Conjecture for large, dense
graphs. Preprint 2008}

\bibitem{diestelBook05}{
R.~Diestel, {\em Graph Theory} (3rd edition). Springer-Verlag, 2005.
}


\bibitem{discrepency}{
 P.~Erd\H{o}s, Z.~F\"uredi, M.~Loebl and V.~T.~S\'os,
Discrepancy of trees. {\em Studia Sci.~Math.~Hungar.}, 30(1-2):47-57, 1995.
}


\bibitem{treeRamsey}
{P.~Haxell, T.~\L uczak and P.~Tingley, Ramsey numbers for trees of small
maximum degree. {\em Combinatorica}, 22(2):287-320, 2002. }


\bibitem{DiaHon}
{J.~Hladk\'y and D.~Piguet,
Loebl--Koml\'os--S\'os conjecture: dense case.
Preprint 2008.
}

\bibitem{komlosSimonovits}
{J.~Koml\'os and M.~Simonovits, Szemer\'edi's regularity lemma and its
appli\-cations in graph theory. In {\em Paul Erd\H os is eighty, (Keszthely,
1993),} volume~2 of {\em Bolyai Soc. Math. Stud.}, pages 295-352, Budapest,
1996. J\'anos Bolyai Math.\ Soc.
}

\bibitem{diam5}{D.~Piguet and M.~Stein, The Loebl--Koml\'os--S\'os conjecture
for trees of diameter $5$ and other special cases. {\em Electr.~J.~of Comb.}, 15:R106, 2008. }


\bibitem{sacwoz}
{J.-F.~Sacl\'e and M.~Wo\'zniak, A note on the Erd\H os--S\'os conjecture for
graphs without ${C}_4$. {\em J.~Combin.~Theory B}, 70(2):229-234, 1997.
}

\bibitem{sze}{
E.~Szemer\'edi, Regular partitions of graphs. {\em Colloques Internationaux
C.N.R.S. 260 -- Probl\`emes Combinatoires et Th\'eorie des Graphes, Orsay},
pages 399-401, 1976.
}


\bibitem{woz}{
M.~Wo\'zniak, On the Erd\H os--S\'os conjecture. {\em J.~Graph Theory},
21(2):229-234, 1996.
}

\bibitem{zhao}
{Y.~Zhao, Proof of the $(n/2-n/2-n/2)$ conjecture for large $n$. Preprint.
}
\end{document}